\documentclass{amsart}

\usepackage{amssymb}

\allowdisplaybreaks

\newtheorem{theorem}{Theorem}[section]
\newtheorem{lemma}[theorem]{Lemma}

\newtheorem{corollary}[theorem]{Corollary}
\newtheorem{remark}[theorem]{Remark}

\theoremstyle{definition}

\newcommand{\dem}{\noindent {\bf Proof. }}
\newcommand{\fin}{\hspace*{\fill} $\square$ \vskip0.2cm}

\newcommand{\Z}{\mathbb{Z}}
\newcommand{\C}{\mathbb{C}}
\newcommand{\X}{\mathrm{X}}
\newcommand{\V}{\mathrm{V}}
\newcommand{\G}{\mathrm{G}}
\newcommand{\Su}{\mathrm{S}}
\newcommand{\Di}{\mathrm{D}}
\newcommand{\OH}{\mathrm{L}}
\newcommand{\Lap}{\mathcal{L}}
\newcommand{\Pol}{\mathbb{P}}
\newcommand{\Rad}{\mathcal{R}_{\subset}}
\newcommand{\Lam}{\Lambda_s}
\newcommand{\lam}{\lambda_s}

\begin{document}

\title[Laplacian operators on Grassmann
graphs] {Laplacian operators and Radon transforms \\ on Grassmann
graphs}

\author[Marco and Parcet]
{Jos\'{e} Manuel Marco and Javier Parcet}

\address{Department of Mathematics, Universidad Aut\'{o}noma de
Madrid, Madrid 28049, Spain}

\email{javier.parcet@uam.es}

\footnote{Partially supported by the Project MTM2004-00678,
Spain.} 

\begin{abstract}
Let $\Omega$ be a vector space over a finite field with $q$
elements. Let $\G$ denote the general linear group of
endomorphisms of $\Omega$ and let us consider the left regular
representation $\rho: \G \rightarrow \mathcal{B}(L_2(\X))$
associated to the natural action of $\G$ on the set $\X$ of linear
subspaces of $\Omega$. In this paper we study a natural basis
$\mathbf{B}$ of the algebra $\mathrm{End_{\G}}(L_2(\X))$ of
intertwining maps on $L_2(\X)$. By using a Laplacian operator on
Grassmann graphs, we identify the kernels in $\mathbf{B}$ as
solutions of a basic hypergeometric difference equation. This
provides two expressions for these kernels. One in terms of the
$q$-Hahn polynomials and the other by means of a Rodrigues type
formula. Finally, we obtain a useful product formula for the
mappings in $\mathbf{B}$. We give two different proofs. One uses
the theory of classical hypergeometric polynomials and the other
is supported by a characterization of spherical functions in
finite symmetric spaces. Both proofs require the use of certain
associated Radon transforms.
\end{abstract}

\maketitle

\section*{Introduction}

Let $\X$ be a homogeneous space with respect to a given finite
group $\G$. That is, the group $\G$ acts transitively on the set
$\X$. Then, we can consider the left regular representation $\rho:
\G \rightarrow \mathcal{B}(\V_{\X})$ of $\G$ into the Hilbert
space $\V_{\X}$ of complex-valued functions $\varphi: \X
\rightarrow \C$. In this context, it is well-known that the
algebra $\mathrm{End_{\G}}(\V_{\X})$ of intertwining operators
with respect to $\rho$ codify some relevant information. For
instance, $\mathrm{End_{\G}}(\V_{\X})$ is abelian if and only if
the left regular representation $\rho$ is multiplicity-free. In
this case, following Terras' book \cite{Te}, we say that $\X$ is a
finite symmetric space with respect to $\G$. When dealing with
finite symmetric spaces, any explicit expression for the kernels
of the orthogonal projections onto the irreducible components of
$\V_{\X}$ is interesting. Indeed, the main motivation lies in the
fact that these expressions can be usually regarded as
combinatorial versions of the irreducible characters. More
generally, if $\X$ and $\mathrm{Y}$ are finite symmetric spaces
with respect to $\G$, some information about the relations between
$\X$ and $\mathrm{Y}$ in the group theory level can be obtained by
studying the space $\mathrm{Hom_{\G}}(\V_{\X},\V_{\mathrm{Y}})$ of
intertwining homomorphisms from $\V_{\X}$ to $\V_{\mathrm{Y}}$.
For instance, the dimension of the space
$\mathrm{Hom_{\G}}(\V_{\X},\V_{\mathrm{Y}})$ coincides with the
number of irreducible components that $\V_{\X}$ and
$\V_{\mathrm{Y}}$ have in common.

\vskip3pt

Two natural problems arise in this setting. In order to state
them, let us consider the Radon transform
$\mathcal{R}_{\null_\mathcal{Z}}: \V_{\X} \rightarrow
\V_{\mathrm{Y}}$ associated to a $\G$-invariant subset
$\mathcal{Z} \subset \X \times \mathrm{Y}$
$$\mathcal{R}_{\null_\mathcal{Z}} \varphi(y) = \sum_{x: \, (x,y)
\in \mathcal{Z}} \varphi(x).$$ Taking the $\G$-invariant subset
$\mathcal{Z}$ to be each of the orbits $\mathcal{O}_1,
\mathcal{O}_2, \ldots, \mathcal{O}_d$ of the action of $\G$ on the
product $\X \times \mathrm{Y}$, we obtain a basis of the space
$\mathrm{Hom_{\G}}(\V_{\X}, \V_{\mathrm{Y}})$ made up of Radon
transforms $$\mathbf{B}_1 = \Big\{ \mathcal{R}_1, \mathcal{R}_2,
\ldots, \mathcal{R}_d \Big\}.$$ On the other hand, since the
dimension of $\mathrm{Hom_{\G}}(\V_{\X},\V_{\mathrm{Y}})$
coincides of the number of irreducible components that $\V_{\X}$
and $\V_{\mathrm{Y}}$ have in common, we obtain the following
orthogonal decompositions
\begin{eqnarray*}
\V_{\X} & = & \mathrm{W}_{\X} \oplus \bigoplus_{1 \le s \le d}
\V_{\X,s}, \\ \V_{\mathrm{Y}} & = & \mathrm{W}_{\mathrm{Y}} \oplus
\bigoplus_{1 \le s \le d} \V_{\mathrm{Y},s},
\end{eqnarray*}
where $\V_{\X,s_1}$ is equivalent to $\V_{\mathrm{Y},s_2}$ if and
only if $s_1 = s_2$. Then we introduce non-zero operators
$\Lambda_s \in \mathrm{Hom_{\G}}(\V_{\X,s}, \V_{\mathrm{Y},s})$
and we regard them as elements of $\mathrm{Hom_{\G}}(\V_{\X},
\V_{\mathrm{Y}})$ vanishing on $\V_{\X} \ominus \V_{\X,s}$. Each
mapping $\Lambda_s$ is an intertwining isomorphism between
$\V_{\X,s}$ and $\V_{\mathrm{Y},s}$ and, by Schur lemma, it is
unique up to a constant factor. This family of mappings provide
another basis of $\mathrm{Hom_{\G}}(\V_{\X}, \V_{\mathrm{Y}})$
$$\mathbf{B}_2 = \Big\{ \Lambda_1, \Lambda_2, \ldots, \Lambda_d
\Big\}.$$

\vskip3pt

Both bases $\mathbf{B}_1$ and $\mathbf{B}_2$ are orthogonal with
respect to the Hilbert-Schmidt inner product. The first problem we
are interested on is to obtain the coefficients relating the bases
described above. On the other hand, given three finite symmetric
spaces $\X_1, \X_2$ and $\X_3$, operator composition provides a
bilinear mapping $$\mathrm{Hom_{\G}}(\V_{\X_2}, \V_{\X_3}) \times
\mathrm{Hom_{\G}}(\V_{\X_1}, \V_{\X_2}) \longrightarrow
\mathrm{Hom_{\G}}(\V_{\X_1}, \V_{\X_3})$$ given by
$$(\Lambda^{2,3}, \Lambda^{1,2}) \mapsto \Lambda^{2,3} \circ
\Lambda^{1,2}.$$ The second problem we want to study is to write
the products $\Lambda_s^{2,3} \circ \Lambda_s^{1,2}$ in terms of
the operators $\Lambda_s^{1,3}$. In this paper we solve the
problems presented above for the Grassmann graphs associated to
the general linear group over a finite field. For the first
problem, we meet a large family of $q$-Hahn polynomials and we
obtain in this way a combinatorial interpretation of this family
of classical hypergeometric polynomials. Besides, we provide a
Rodrigues type formula for these polynomials adapted to the
present framework. For the second problem, we obtain the product
formula by two different processes. One lies in the theory of
classical hypergeometric polynomials and the other in the theory
of finite symmetric spaces and spherical functions. The second
approach also provides certain identities for the kernels in
$\mathbf{B}_2$ which might be of independent interest. The results
we present in this paper constitute the $q$-analogue of those
given in \cite{MP1}. The main idea is to identify certain
Laplacian type operators on Grassmann graphs as hypergeometric
type operators. This procedure will allow us to apply some results
on classical hypergeometric polynomials which have appeared
recently in \cite{MP2}. The paper \cite{MP2} provides a new
approach to the theory of classical hypergeometric polynomials
which somehow lives between the theories developed by Askey and
Wilson \cite{AW} on one side and by Nikiforov, Suslov and Uvarov
\cite{NSU} on the other. One of the main motivations for this
paper is to show the efficiency of the point of view suggested in
\cite{MP2}.

\vskip3pt

The organization of the paper is as follows. Let $\Omega$ be a
finite-dimensional vector space over a finite field $\mathbb{K}$
and let $\G$ be the general linear group of endomorphisms of
$\Omega$. Let us consider the left regular representation $\rho:
\G \rightarrow \mathcal{B}(\V_{\X})$ associated with the natural
action of $\G$ on the set $\X$ of linear subspaces of $\Omega$.
Section \ref{Section1} is devoted to describe a natural basis
$\mathbf{B}$ of the intertwining algebra
$\mathrm{End_{\G}}(\V_{\X})$. The Laplacian operators on Grassmann
graphs are studied in Section \ref{Section2}. This is used in
Section \ref{Section3} to identify these operators with some
well-known hypergeometric type operators. Then we provide
polynomic expressions and Rodrigues type formulas for the kernels
of the operators in $\mathbf{B}$. Finally, in Sections
\ref{Section4} and \ref{Section5} we give two different proofs of
the product formula mentioned above. Section \ref{Section4} uses
the theory of classical hypergeometric polynomials while the proof
given in Section \ref{Section5} lies in a characterization of
spherical functions on finite symmetric spaces.

\vskip3pt

After having written this paper, the authors were informed on the
existence of Dunkl's paper \cite{D}, which also identifies the
kernels mentioned above as $q$-Hahn polynomials and studies
similar relations for them. We note however that there also exist
some significant differences between both papers. Indeed, our main
product formula is not obtained in \cite{D} while the Rodrigues
formulas deduced from \cite{MP2} (with non-ramified weights) are
new.

\section{The object to study}
\label{Section1}

Let $\mathbb{K}$ denote the field $\mathbb{F}_q$ with $q$ elements
for some power of a prime $q$. In what follows, $\Omega$ will be a
finite-dimensional vector space over $\mathbb{K}$. If $n$ stands
for $\dim \Omega$, we shall also consider the set $\X$ of linear
subspaces of $\Omega$ and the sets $\X_r$ of $r$-dimensional
subspaces of $\Omega$ for $0 \le r \le n$. For any $(x,y) \in \X
\times \X$, we define
\begin{eqnarray*}
\partial(x,y) & = & \dim \left( x/(x \cap y) \right), \\
\overline{\partial}(x,y) & = & \partial(x,y) + \partial(y,x).
\end{eqnarray*}
Recall that, by Grassmann formula, we also have $\partial(x,y) =
\dim((x+y)/y)$. Taking $\partial(x) = \dim (x)$, we can write
$\overline{\partial}(x,y) = \partial(x) + \partial(y) - 2
\partial(x \cap y)$. The function $\overline{\partial}: \X \times
\X \rightarrow \mathbb{R}_+$ is clearly a graph distance on $\X$
and the same happens with the restriction $\partial_r: \X_r \times
\X_r \rightarrow \mathbb{R}_+$ of $\partial$ on each $\X_r$. The
distance $\partial_r$ imposes on $\X_r$ an structure of
distance-regular graph. These graphs are well-known in the
literature as \textbf{Grassmann graphs}, see \cite{BCN} for more
on this. Besides, we shall also need to consider the vector space
$\V$ of complex valued functions $\varphi: \X \rightarrow \C$ and
its subspaces $\V_r$ made up of functions $\varphi: \X_r
\rightarrow \C$ for $0 \le r \le n$. Notice that both $\V$ and
$\V_r$ are vector spaces over the complex field. If we consider
the natural Hilbert space structure on these spaces, so that $\V =
L_2(\X)$ and $\V_r = L_2(\X_r)$, we clearly have the orthogonal
decomposition $$\V = \bigoplus_{r=0}^n \V_r.$$

\subsection{Finite symmetric spaces}
\label{Subsection1.1}

Before stating in detail the problem we want to study, we give a
brief summary of results on finite symmetric spaces and spherical
functions that will be used in the sequel. For further information
on these topics see \cite{Te} and the references cited there. Let
$\G$ be a finite group acting on a finite set $\X$. This action
gives rise to the left regular representation $\rho: \G
\rightarrow \mathcal{B}(L_2(\X))$, defined as follows $$\big(
\rho(g) \varphi \big) (x) = \varphi (g^{-1} x).$$ Assume that the
action $\G \times \X \rightarrow \X$ is transitive, so that $\X$
becomes a homogeneous space. Then $\X$ is called a \textbf{finite
symmetric space} with respect to the group $\G$ if the algebra
$\mathrm{End_{\G}}(L_2(\X))$ of intertwining endomorphisms of
$L_2(\X)$ is abelian.

\begin{remark} \label{Remark-Multiplicity-Free}
\emph{$\mathrm{End_{\G}}(L_2(\X))$ is abelian if and only if
$\rho$ is multiplicity-free, see \cite{Te}.}
\end{remark}

Now, if we are given two finite symmetric spaces $\X_1$ and $\X_2$
with respect to $\G$, let us denote by $\rho_1$ and $\rho_2$ the
corresponding associated unitary representations. Then we assign
to each $\Lambda$ in $\mathrm{Hom} (L_2(\X_1),L_2(\X_2))$ its
kernel $\lambda: \X_2 \times \X_1 \rightarrow \C$ with respect to
the natural bases. The mapping $\Lambda \rightarrow \lambda$ is
clearly a linear isomorphism $\Psi: \mbox{Hom}
(L_2(\X_1),L_2(\X_2)) \rightarrow L_2(\X_2 \times \X_1)$ with
$\Lambda$ and $\lambda$ related by $$\Lambda \varphi(x_2) =
\sum_{x_1 \in \X_1}^{\null} \lambda(x_2,x_1) \varphi(x_1).$$ Let
$\mathrm{Hom_{\G}}(L_2(\X_1),L_2(\X_2))$ be the space of
intertwining maps for $\rho_1$ and $\rho_2$. If we compare
$\Lambda \circ \rho_1(g)$ and $\rho_2(g) \circ \Lambda$ written in
terms of $\lambda$, it is not difficult to check that $\Lambda \in
\mathrm{Hom_{\G}}(L_2(\X_1),L_2(\X_2))$ if and only if
$\lambda(gx_2,gx_1) = \lambda(x_2,x_1)$ holds for all $g \in \G$
and all $(x_1,x_2) \in \X_1 \times \X_2$. That is, $\Lambda$ is an
intertwining operator for $\rho_1$ and $\rho_2$ if and only if
$\lambda$ is constant on the orbits of the action $$\G \times \X_2
\times \X_1 \ni (g, (x_2,x_1)) \longmapsto (g x_2, g x_1) \in \X_2
\times \X_1.$$

Now, assume we are given a transitive action $\G \times \X
\rightarrow \X$ of a finite group $\G$ on a finite set $\X$
endowed with a distance $\partial$. We say that $\X$ is a
\textbf{two-point homogeneous space} when for any two pairs
$(x_1,x_2),(y_1,y_2) \in \X \times \X$ satisfying
$\partial(x_1,x_2) = \partial(y_1,y_2)$, there exists $g \in \G$
such that $g x_1 = y_1$ and $g x_2 = y_2$.

\begin{remark} \label{Remark-Symmetric-Action}
\emph{An action $\G \times \X \rightarrow \X$ is called
\textbf{symmetric} when for any $x_1,x_2 \in \X$ there exists $g
\in \G$ such that $g x_1 = x_2$ and $g x_2 = x_1$. Any two-point
homogeneous space $\X$ is clearly equipped with a symmetric action
and in that case $\X$ becomes a finite symmetric space with
respect to $\G$. Namely, if the action of $\G$ on $\X$ is
symmetric, then $\Psi \left( \mathrm{End_{\G}}(L_2(\X)) \right)$
is a subalgebra of $L_2(\X \times \X)$ made up of symmetric
matrices, hence abelian. Finally, since the mapping $\Psi$ is an
algebra isomorphism when $\X_1 = \X_2$, it turns out that $\X$ is
a finite symmetric space.}
\end{remark}

Let us write $\widehat{\G}$ for the dual object of $\G$. That is,
the set of irreducible unitary representations of $\G$ up to
unitary equivalence. Let us consider the set $$\widehat{\G}_{\X} =
\Big\{ \pi \in \widehat{\G}: \, \mathrm{Mult}_{\pi} (\rho) \neq 0
\Big\}.$$ Notice that if $\X$ is symmetric with respect to $\G$,
then every $\pi \in \widehat{\G}_{\X}$ satisfies
$\mbox{Mult}_{\pi} (\rho) = 1$ since $\rho$ is multiplicity-free
by Remark \ref{Remark-Multiplicity-Free}. This set allows us to
decompose $L_2(\X)$ into irreducible components $$L_2(\X) =
\bigoplus_{\pi \in \widehat{\G}_{\X}} L_2(\X)_{\pi}.$$ We denote
by $\mathrm{P}_{\pi}$ the orthogonal projection onto
$L_2(\X)_{\pi}$. The kernel of $\mathrm{P}_{\pi}$ will be denoted
by $p_{\pi}$. The \textbf{spherical functions} on $\X$ are defined
by $$\psi_{\X,\pi} = \frac{|\X|}{d(\pi)} p_{\pi} \in \Psi \left(
\mathrm{End_{\G}} (L_2(\X)) \right),$$ where $\pi \in
\widehat{\G}_{\X}$ and $d(\pi)$ denotes the degree of $\pi$. We
shall also write $\mathrm{S}_{\X,\pi}$ for the associated operator
in $\mathrm{End_{\G}}(L_2(\X))$ with kernel $\psi_{\X,\pi}$. A
slightly modified version of the following result can be found in
Terra's book \cite[Th. 1 of Chapter 20]{Te}.

\begin{theorem} \label{Theorem-Spherical}
Let $\X$ be a finite symmetric space with respect to the finite
group $\G$ and let $\psi \in \Psi \left( \mathrm{End_{\G}}
(L_2(\X)) \right)$, then the following are equivalent:
\begin{itemize}
\item[$\mathrm{(a)}$] There exists $\pi \in \widehat{\G}_{\X}$ such
that $\psi = \psi_{\X,\pi}$.
\item[$\mathrm{(b)}$] The function $\psi$ satisfies
$\psi(x_0,x_0) = 1$ for all $x_0 \in \X$ and
$$\frac{1}{|\G_{x_0}|} \sum_{g \in \G_{x_0}} \psi(gx_1,x_2) =
\psi(x_1,x_0) \psi(x_0,x_2) = \frac{1}{|\G_{x_0}|} \sum_{g \in
\G_{x_0}} \psi(x_1,gx_2)$$ for every $x_1,x_2 \in \X$ and where
$\G_{x_0}$ denotes the isotropy subgroup of $x_0$.
\end{itemize}
\end{theorem}

\subsection{Notation and results from $q$-combinatorics}
\label{Subsection1.2}

We shall also need some results from $q$-combinatorics that we
summarize here. Our notation will follow the book \cite{GR} by
Gasper and Rahman. For some related results in $q$-combinatorics,
the reader is referred to \cite{VW}. The \textbf{$q$-shifted
factorials} are defined as follows $$(u;q)_k = \prod_{j=0}^{k-1}
(1-q^j u) \qquad \mbox{and} \qquad (u;q^{-1})_k =
\prod_{j=0}^{k-1} (1-q^{-j} u).$$ Then, the \textbf{$q$-binomial
coefficients} can be written as $$\Big[ \!\!
\begin{array}{c} m \\ k \end{array} \!\! \Big]_q =
\frac{(q;q)_m}{(q;q)_k (q;q)_{m-k}} =
\frac{(q^m;q^{-1})_k}{(q^k;q^{-1})_k}= \Big[ \!\!
\begin{array}{c} m \\ m-k \end{array} \!\! \Big]_q.$$ Now we
present some well-known combinatorial identities that will be used
in this paper with no further comment. Let us recall some of the
objects introduced above: $\mathbb{K}$, $\Omega$, $\X_r$,
$\mathrm{V}_r$, $\partial, \ldots$ Besides, let us consider the
general linear group $\mathrm{GL}(n,\mathbb{K})$. That is, the
group of endomorphisms of $\Omega$. Then, we have
\begin{eqnarray}
\label{Equation-Cardinal-Group} |\mathrm{GL}(n,\mathbb{K})| & = &
(-1)^n q^{{{n} \choose {2}}} (q;q)_n, \\
\label{Equation-Cardinal-Graph} \dim \V_r & = & \Big[ \!\!
\begin{array}{c} n \\ r \end{array} \!\! \Big]_q.
\end{eqnarray}
Moreover, given $0 \le r_1 \le r \le r_2 \le n$ and $(x_1,x_2) \in
\X_{r_1} \times \X_{r_2}$, we have
\begin{equation} \label{Equation-Cardinal-3}
\Big| \Big\{ x \in \X_r \, \big| \ x_1 \subset x \subset x_2
\Big\} \Big| = \Big[ \!\! \begin{array}{c} r_2 - r_1 \\ r - r_1
\end{array} \!\! \Big]_q.
\end{equation}
Eventually, we shall also use the \textbf{$q$-multinomial
coefficients} $$\Big[ \!\! \begin{array}{c} m \\ r_1, r_2, \ldots,
r_k \end{array} \!\! \Big]_q = \frac{(q;q)_m}{(q;q)_{r_1} \cdots
(q;q)_{r_k}}$$ with $m = \sum_1^k r_j$ and $r_1, r_2, \ldots, r_k
\ge 0$. The combinatorial interpretation is $$\Big| \Big\{ (x_1,
\ldots, x_k) \in \prod_{j=1}^k \X_{r_j} \, \big| \ \Omega =
\bigoplus_{j=1}^k x_j \Big\} \Big| = \Big[ \!\!
\begin{array}{c} n \\ r_1, \ldots, r_k \end{array} \!\! \Big]_q
\prod_{i<j} q^{r_ir_j}.$$ In particular, given $z \in \X_r$ we
have
\begin{equation} \label{Equation-Suplement}
\Big| \Big\{w \in \X_{n-r} \, \big| \ \Omega = z \oplus w \Big\}
\Big| = q^{r(n-r)}.
\end{equation}
Let $\mathrm{I}_n(r_1,r_2)$ be the set of all possible values of
the parameter $\partial(x_2,x_1)$ for $x_1 \in \X_{r_1}$ and $x_2
\in \X_{r_2}$. Given $t \in \mathrm{I}_n(r_1,r_2)$, the following
identity follows from the relations above and will be very useful
for our forthcoming computations
\begin{eqnarray} \label{Equation-Distances}
\lefteqn{\Big| \Big\{ (x_1,x_2) \in \X_{r_1} \times \X_{r_2} \,
\big| \ \partial(x_2,x_1) = t \Big\} \Big|} \\ \nonumber  & = &
q^{t(r_1-r_2+t)} \Big[ \!\! \begin{array}{c} n \\ t, r_2-t, r_1 -
r_2 + t, n - r_1 - t \end{array} \!\! \Big]_q.
\end{eqnarray}

\begin{remark}
\emph{We shall also use the notation $\displaystyle (u_1, u_2,
\ldots, u_d; q)_k = \prod_{j=1}^d (u_j;q)_k$.}
\end{remark}

\subsection{The basis of the algebra
$\mathrm{End_{\G}}(\V)$} \label{Subsection1.3}

Let $\G$ be the general linear group $\mathrm{GL}(n,\mathbb{K})$
of endomorphisms of $\Omega$ considered above. This group acts
naturally on the set $\X$ of linear subspaces of the vector space
$\Omega$. The orbits of this action are the subsets $\X_0, \X_1,
\ldots, \X_n$ of linear subspaces of dimensions $0,1, \ldots, n$.

\begin{remark}
\emph{The action of $\G$ on $\X_r$ preserves $\partial_r$ and
$\X_r$ is clearly a two-point homogeneous space for each $0 \le r
\le n$. In particular, by Remark \ref{Remark-Symmetric-Action},
$\X_r$ becomes a finite symmetric space with respect to $\G$.
Besides, Remark \ref{Remark-Multiplicity-Free} gives that the left
regular representations of $\G$ into $\V_r$ are multiplicity-free
for each $0 \le r \le n$.}
\end{remark}

Given any two integers $0 \le r_1, r_2 \le n$, we shall identify
each space $\mathrm{Hom}(\V_{r_1},\V_{r_2})$ with a subspace of
$\mathrm{End}(\V)$ by right multiplication by the orthogonal
projection from $\V$ onto $\V_{r_1}$. Applying the same
identification for the intertwining operators, we obtain the
following decompositions
\begin{eqnarray*}
\mathrm{End}(\V) & = & \bigoplus_{r_1,r_2}
\mathrm{Hom}(\V_{r_1},\V_{r_2}), \\ \mathrm{End_{\G}}(\V) & = &
\bigoplus_{r_1,r_2} \mathrm{Hom_{\G}}(\V_{r_1},\V_{r_2}).
\end{eqnarray*}
The algebra $\mathrm{End}(\V)$ is a Hilbert space with respect to
the Hilbert-Schmidt inner product and the direct sums given above
become orthogonal decompositions with respect to this structure.
Besides, we know that the kernel of any intertwining operator
$\Lambda  \in \mathrm{Hom_{\G}}(\V_{r_1},\V_{r_2})$ is constant on
the orbits of $(g,(x_2,x_1)) \mapsto (gx_2,gx_1)$. These orbits
are completely determined by the parameter $\partial(x_2,x_1)$. In
particular, the kernel of $\Lambda$ can be regarded as a function
$\lambda: \mathrm{I}_n(r_1,r_2) \rightarrow \C$, where the index
set $\mathrm{I}_n(r_1,r_2)$ was considered above
$$\mathrm{I}_n(r_1,r_2) = \Big\{
\partial(x_2,x_1) \, \big| \ x_1 \in \X_{r_1}, x_2 \in \X_{r_2}
\Big\} = \Big\{0 \vee (r_2 - r_1) \le t \le r_2 \wedge (n-r_1)
\Big\}.$$ Here $\wedge$ stands for $\min$ and $\vee$ for $\max$.
Therefore we have $$\Lambda \varphi(x_2) = \sum_{x_1 \in
\X_{r_1}}^{\null} \lambda(\partial(x_2,x_1)) \, \varphi(x_1)$$ for
any $\varphi \in \V_{r_1}$ and $$\dim \big( \mathrm{Hom_{\G}}
(\V_{r_1},\V_{r_2}) \big) = \big| \mathrm{I}_n(r_1,r_2) \big| = 1
+ \mathrm{N}(r_1,r_2),$$ where $\mathrm{N}(r_1,r_2) = r_1 \wedge
r_2 \wedge (n-r_1) \wedge (n-r_2)$. Reciprocally, any $\lambda:
\mathrm{I}_n(r_1,r_2) \rightarrow \C$ determines an operator
$\Lambda \in \mathrm{Hom_{\G}}(\V_{r_1},\V_{r_2})$. Moreover,
since $\V_r$ is multiplicity free, Schur lemma gives that the
dimension of $\mathrm{Hom_{\G}} (\V_{r_1},\V_{r_2})$ is the number
of irreducible components that $\V_{r_1}$ and $\V_{r_2}$ have in
common. In particular, we have
\begin{itemize}
\item $\V_r$ has $1 + \mathrm{N}(r,r) = 1 + r \wedge (n-r)$
irreducible components.
\item The number of irreducible components that $\V_{r_1}$ and
$\V_{r_2}$ have in common is the minimum of the numbers of
irreducible components of $\V_{r_1}$ and $\V_{r_2}$.
\end{itemize}
Therefore, there exist a family of inequivalent irreducible
unitary representations $\pi_s: \G \rightarrow
\mathcal{B}(\mathcal{H}_s)$ such that, if we denote by $\V_{r,s}$
the $\G$-invariant subspace of $\V_r$ equivalent to
$\mathcal{H}_s$, the left regular representation $\rho_r: \G
\rightarrow \mathcal{B}(\V_r)$ and the Hilbert space $\V_r$
decompose into irreducibles as follows $$\rho_r \simeq
\bigoplus_{s=0}^{r \wedge (n-r)} \pi_s \qquad \mbox{and} \qquad
\V_r = \bigoplus_{s=0}^{r \wedge (n-r)} \V_{r,s}.$$ Moreover, the
representations of $\G$ into $\V_{r_1,s_1}$ and $\V_{r_2,s_2}$ are
equivalent if and only if $s_1=s_2$. Finally we note that
\begin{equation} \label{Equation-Dimension}
\dim \V_{r,s} = \dim \V_{s,s} = \dim \V_s - \dim \V_{s-1} = \Big[
\!\! \begin{array}{c} n \\ s \end{array} \!\! \Big]_q - \Big[ \!\!
\begin{array}{c} n \\ s-1 \end{array} \!\! \Big]_q.
\end{equation}
Here we assume by convention $\X_{-1} = \emptyset$, so that $\dim
\V_{-1} = 0$. The last identity in (\ref{Equation-Dimension})
follows from relation (\ref{Equation-Cardinal-Graph}). By Schur
lemma we know that $\mathrm{Hom_{\G}}(\V_{r_1,s},\V_{r_2,s})$ is
one-dimensional. In summary, we have obtained an orthogonal
decomposition of the algebra of intertwining operators
$\mathrm{End_{\G}}(\V)$ into one-dimensional subspaces
$$\mathrm{End_{\G}}(\V) = \bigoplus_{0 \le r_1,r_2 \le n}
\bigoplus_{s=0}^{\mathrm{N}(r_1,r_2)}
\mathrm{Hom_{\G}}(\V_{r_1,s},\V_{r_2,s}).$$ This decomposition
provides a natural basis of the algebra $\mathrm{End_{\G}}(\V)$
which will be the object of our study. Namely, taking a non-zero
element $\Lambda_s^{r_1,r_2}$ in each space
$\mathrm{Hom_{\G}}(\V_{r_1,s}, \V_{r_2,s})$, we obtain a basis
$\mathbf{B}$ of the algebra $\mathrm{End_{\G}}(\V)$. Our
definition of $\Lambda_s^{r_1,r_2}$ is still ambiguous since we
have only defined it up to a constant factor. We shall precise
this below. As it was announced in the Introduction, the aim of
this paper is to provide several expressions for the kernels
$\lambda_s^{r_1,r_2}$ and to give an explicit formula for the
mapping product $\Lambda_s^{r_2,r_3} \circ \Lambda_s^{r_1,r_2}$.

\section{Laplacian operators on graphs}
\label{Section2}

In this section we deal with some Laplacian type operators which
will be useful to identify certain difference equation satisfied
by the kernels of $\Lambda_s^{r_1,r_2}$ in the usual rank of
parameters for $r_1,r_2$ and $s$. We begin by recalling some
general results for Laplacian operators on graphs. Then we focus
on the Grassmann graphs.

\subsection{General results}
\label{Subsection2.1}

Let $\X$ be a finite distance-regular graph and let $\partial$ be
the distance on $\X$. Assume there exists a finite group $\G$
acting on $\X$ such that the graph $\X$ becomes a two-point
homogeneous space with respect to $\G$. This structure on $\X$
allows us to define two Laplacian type operators on the vector
space $\V$ of complex valued functions $\varphi: \X \rightarrow
\C$. First, given $x \in \X$, we consider the set
$$\mathbf{S}_1(x) = \Big\{ y \in \X \, \big| \
\partial(x,y)=1 \Big\}.$$ The regularity for the distance imposed
on the graph $\X$ implies that the cardinality of the set
$\mathbf{S}_1(x)$ does not depend on $x$. This cardinality
$\mathrm{val}(\X)$ is called in the literature the valence of
$\X$. The \textbf{graph Laplacian} $\Lap_{\X}: \V \rightarrow \V$
is the operator defined as follows $$\Lap_{\X} \varphi(x) =
\sum_{y \in \mathbf{S}_1(x)}^{\null} \big( \varphi(y) - \varphi(x)
\big) = \sum_{y \in \mathbf{S}_1(x)} \varphi(y) - \mathrm{val}(\X)
\varphi(x).$$ Second, let us consider a subset $\mathbf{T}$ of
$\G$ satisfying the following properties
\begin{itemize}
\item $\mathbf{T} = \mathbf{T}^{-1}$.
\item $\mathbf{T} = g \mathbf{T} g^{-1}$ for all
$g \in \G$.
\item $\partial(x,hx) \le 1$ for all $x \in \X$
and all $h \in \mathbf{T}$.
\item There exists $x \in \X$ such that $\mathbf{T} \nsubseteq \G_x$,
the isotropy subgroup of $x$.
\end{itemize}
The \textbf{group Laplacian} $\Lap_{\mathbf{T},\X}: \V \rightarrow
\V$ is defined as follows $$\Lap_{\mathbf{T},\X} \varphi(x) =
\sum_{h \in \mathbf{T}}^{\null} \big( \varphi(hx) - \varphi(x)
\big) = \sum_{h \in \mathbf{T}} \varphi(hx) - |\mathbf{T}|
\varphi(x).$$

Both the graph and the group Laplacians are self-adjoint operators
with respect to the natural inner product on $\V$. This is an easy
exercise that we leave to the reader. Besides, recalling that an
endomorphism of $\V$ is an intertwining operator if and only if
its kernel is constant on the orbits of the action of $\G$ on $\X
\times \X$, it is not difficult to check that both the graph and
the group Laplacians belong to the intertwining algebra
$\mathrm{End_{\G}}(\V)$. Now let $\mathbf{T}$ be a subset of $\G$
satisfying the properties above. Then, given $x \in \X$ and $y \in
\mathbf{S}_1(x)$, we define
\begin{eqnarray*}
\mathbf{T}_x & = & \Big\{ h \in \mathbf{T} \, \big| \ hx = x
\Big\}, \\ \mathbf{T}_{x,y} & = & \Big\{ h \in \mathbf{T} \, \big|
\ hx=y \Big\}.
\end{eqnarray*}
From the properties of $\mathbf{T}$ it follows that $g
\mathbf{T}_x g^{-1} = \mathbf{T}_{gx}$ and $g \mathbf{T}_{x,y}
g^{-1} = \mathbf{T}_{gx,gy}$ for all $g \in \G$. In particular,
since the distance-regular graph $\X$ is assumed to be a two-point
homogeneous space, the numbers
\begin{eqnarray*}
\gamma_0(\mathbf{T},\X) & = & |\mathbf{T}_x|, \\
\gamma_1(\mathbf{T},\X) & = & |\mathbf{T}_{x,y}|,
\end{eqnarray*}
do not depend on the election of $x \in \X$ and $y \in
\mathbf{S}_1(x)$. We now state some basic results on these
Laplacian operators that will be used in the sequel. Notice that
the fourth condition imposed on $\mathbf{T}$ implies that
$\gamma_1(\mathbf{T},\X) > 0$.

\begin{lemma} \label{Lemma-Laplacian-Ratio}
$|\mathbf{T}| = \gamma_0(\mathbf{T},\X) + \mathrm{val}(\X) \,
\gamma_1(\mathbf{T},\X)$ and $\Lap_{\mathbf{T},\X} =
\gamma_1(\mathbf{T},\X) \, \Lap_{\X}$.
\end{lemma}

\dem Obviously we have $$\mathbf{T} = \mathbf{T}_x \cup \Big[
\bigcup_{y \in \mathbf{S}_1(x)} \mathbf{T}_{x,y} \Big]$$ with
disjoint unions. Hence, the first assertion follows. Besides, we
notice that $$\Lap_{\mathbf{T},\X} \varphi(x) = \sum_{h \in
\mathbf{T}} \big( \varphi(hx) - \varphi(x) \big) = \sum_{y \in
\mathbf{S}_1(x)} \sum_{h \in \mathbf{T}_{x,y}} \big( \varphi(hx) -
\varphi(x) \big).$$ By definition, the last expression is
$\gamma_1(\mathbf{T},\X) \, \Lap_{\X} \varphi(x)$. This concludes
the proof. \fin

\begin{lemma} \label{Lemma-Laplacian-Irreducible}
Let $\mathrm{W}$ be a $\G$-invariant subspace of $\V$. Then, the
group Laplacian $\Lap_{\mathbf{T},\X}$ preserves $\mathrm{W}$.
Besides, if $\mathrm{W}$ is irreducible, there exists a complex
number $\mu$ depending only on the representation of $\G$ into
$\mathrm{W}$ such that $${\Lap_{\mathbf{T},\X}}_{|_{\mathrm{W}}} +
\mu \, 1_{\mathrm{W}} = 0.$$
\end{lemma}

\dem Let $\rho: \G \rightarrow \mathcal{B}(\V)$ be the left
regular representation associated to the action of $\G$ into $\X$.
Then, for any $\varphi \in \mathrm{W}$, we can write
$$\Lap_{\mathbf{T},\X} \varphi(x) = \sum_{h \in \mathbf{T}} \big(
\rho(h) \varphi(x) - 1_{\V} \varphi(x) \big).$$ The first claim
follows from the relation above. The second is a consequence of
Schur lemma. Namely, if $\mathrm{W}$ is irreducible and $\pi: \G
\rightarrow \mathcal{B}(\mathrm{W})$ denotes the restriction of
$\rho$ to $\mathrm{W}$, we have
$${\Lap_{\mathbf{T},\X}}_{|\mathrm{W}} = \sum_{h \in \mathbf{T}}
\big( \pi(h) - 1_{\mathrm{W}} \big).$$ Therefore, since this is an
intertwining operator with respect to $\pi$, Schur lemma gives
that ${\Lap_{\mathbf{T},\X}}_{|\mathrm{W}} + \mu \, 1_{\mathrm{W}}
= 0$ with $\mu$ depending only on the representation $\pi$. \fin

\subsection{Laplacian operators on the graphs $\mathrm{X}_r$}
\label{subsection2.2}

In this paragraph we return to the study of the intertwining
algebra $\mathrm{End_{\G}}(\V)$ described in Section
\ref{Section1}. In particular, we shall work with the general
linear group $\G$ of endomorphisms of $\Omega$ and the Grassmann
graphs $\X_r$ for $0 \le r \le n$. We denote by $\Lap_r: \V_r
\rightarrow \V_r$ the graph Laplacian on $\X_r$ while
$\Lap_{\mathbf{T},r}: \V_r \rightarrow \V_r$ stands for the group
Laplacian. If $\mathrm{rk}(\Lambda)$ denotes the rank of a mapping
$\Lambda$, we consider the subset $\mathbf{T}$ of $\G$ defined  as
follows
\begin{eqnarray*}
\mathbf{T} & = & \Big\{ g \in \G \, \big| \ \mathrm{rk}(g-1) = 1,
(g-1)^2 = 0 \Big\} \\ & = & \Big\{ 1 + \omega \otimes \alpha \,
\big| \ \alpha \in \Omega^* \setminus \{0\}, \omega \in
\mathrm{Ker}(\alpha) \setminus \{0\} \Big\},
\end{eqnarray*}
where $\omega \otimes \alpha \in \mathrm{End}(\Omega)$ is given by
$(\omega \otimes \alpha)(\omega_0) = \alpha(\omega_0) \omega$. We
need to check that the subset $\mathbf{T}$ satisfies the
properties introduced before the definition of the group
Laplacian. But this is an easy exercise that we leave to the
reader. The main results of this section are summarized in the
following theorem. We shall also need to consider the operator
$\Lap_{r_1,r_2}: \mathrm{Hom_{\G}}(\V_{r_1},\V_{r_2}) \rightarrow
\mathrm{Hom_{\G}}(\V_{r_1},\V_{r_2})$ defined by the following
relation $$\Lap_{r_1,r_2} \Lambda = \Lap_{r_2} \circ \Lambda.$$

\begin{theorem} \label{Theorem-Laplacian-Grassmann}
The Laplacian operators considered above satisfy:
\begin{itemize}
\item[$\mathrm{(a)}$] We have $$\gamma_1(\mathbf{T},\X_r) =
q^{n-2} (q-1) \qquad \mbox{and} \qquad \Lap_{\mathbf{T},r} =
q^{n-2} (q-1) \Lap_r.$$
\item[$\mathrm{(b)}$] The mappings $\Lambda \in
\mathrm{Hom_{\G}}(\V_{r_1},\V_{r_2})$ satisfy $$\Lap_{r_2} \circ
\Lambda = \Lap_{r_1,r_2} \Lambda = \Lambda \circ \Lap_{r_1}.$$
\item[$\mathrm{(c)}$] Given $\varphi \in \V_{r,s}$ and $\Lambda \in
\mathrm{Hom_{\G}}(\V_{r_1,s},\V_{r_2,s})$, there exists $\mu_s(n)$
such that
\begin{eqnarray*}
\Lap_r \varphi(x) + \mu_s(n) \varphi(x) & = & 0, \\ \Lap_{r_1,r_2}
\Lambda + \mu_s(n) \Lambda & = & 0.
\end{eqnarray*}
\item[$\mathrm{(d)}$] If $\Lambda \in
\mathrm{Hom_{\G}}(\V_{r_1},\V_{r_2})$ has kernel $\lambda:
\mathrm{I}_n(r_1,r_2) \rightarrow \C$, then $$\Lambda' =
\Lap_{r_1,r_2} \Lambda$$ has kernel $\lambda':
\mathrm{I}_n(r_1,r_2) \rightarrow \C$ given by $$\lambda'(t) =
\mathbf{b}_{r_1,r_2}(t) \big( \lambda(t+1) - \lambda(t) \big) +
\mathbf{c}_{r_1,r_2}(t) \big( \lambda(t-1) - \lambda(t) \big),$$
where, given $(x_1,x_2) \in \X_{r_1} \times \X_{r_2}$ with
$\partial(x_2,x_1) = t$, we have
\begin{eqnarray*}
\mathbf{b}_{r_1,r_2}(t) & = & \Big| \Big\{ y \in \X_{r_2} \, \big|
\ \partial(y,x_2) = 1, \partial(y,x_1) = t+1 \Big\} \Big|,
\\ \mathbf{c}_{r_1,r_2}(t) & = & \Big| \Big\{ y \in \X_{r_2} \,
\big| \ \partial(y,x_2) = 1, \partial(y,x_1) = t-1 \Big\} \Big|.
\end{eqnarray*}
\item[$\mathrm{(e)}$] The following expressions hold
\begin{eqnarray*}
\mathbf{b}_{r_1,r_2}(t) & = & \frac{q^{r_1-r_2+1}
(q^t-q^{r_2})(q^t-q^{n-r_1})}{(q-1)^2}, \\ \mathbf{c}_{r_1,r_2}(t)
& = & \frac{q^{r_1-r_2} (q^t-1)(q^t-q^{r_2-r_1})}{(q-1)^2}.
\end{eqnarray*}
\end{itemize}
\end{theorem}

The proof of Theorem \ref{Theorem-Laplacian-Grassmann} requires
several auxiliary results. We shall state and prove these results
as they are needed. Given $x \in \X_r$, we consider the sets
\begin{eqnarray*}
\mathbf{S}_1(x) & = & \Big\{y \in \X_r \, \big| \ \partial_r(x,y)
= 1 \Big \}, \\ \widetilde{\mathbf{S}}_1(x) & = & \Big\{
(z,\omega) \in \X_{r-1} \times \Omega \, \big| \ z \subset x,
\omega \notin x \Big\}.
\end{eqnarray*}

\begin{lemma} \label{Lemma-Covering}
Given any $x \in \X_r$, the mapping $$(z,\omega) \in
\widetilde{\mathbf{S}}_1(x) \longmapsto y = z \oplus \mathbb{K}
\omega \in \mathbf{S}_1(x)$$ is surjective, $z = x \cap y$ and
there exist $q^{r-1}(q-1)$ possible $\omega$ for each $y \in
\mathbf{S}_1(x)$.
\end{lemma}

\dem Given $x,y \in \X_r$ we have $\partial_r(x,y) = r -
\partial(x \cap y)$. Hence $\partial_r(x,y) = 1$ is equivalent to
$\partial(x \cap y) = r - 1$. In other words, $y = z \oplus
\mathbb{K} \omega$ with $$z = x \cap y \qquad \mbox{and} \qquad
\omega \in y \setminus z.$$ Hence we have $\omega \notin x$,
$(z,\omega) \in \widetilde{\mathbf{S}}_1(x)$ and there exists $$|y
\setminus z| = q^r - q^{r-1}$$ possible elections for the vector
$\omega$. Therefore, the proof is completed. \fin

Lemma \ref{Lemma-Covering} allows us to compute the valence  of
the Grassmann graphs $\X_r$. Namely, since $|\X_r| = \dim
\mathrm{V}_r$, we can apply (\ref{Equation-Cardinal-Graph}) to
obtain  for $x \in \X_r$
\begin {equation} \label{Equation-Valence}
\mathrm{val}(\X_r) = |\mathbf{S}_1(x)| =
\frac{|\widetilde{\mathbf{S}}_1(x)|}{q^{r-1}(q-1)} = \Big[ \!\!
\begin{array}{c} r \\ 1 \end{array} \!\! \Big]_q
\frac{q^n - q^r}{q^{r-1}(q-1)} = \frac{q (q^r - 1) (q^{n-r}
-1)}{(q-1)^2}.
\end{equation}

\begin{lemma} \label{Lemma-Isotropy}
Given $g = 1 + \omega \otimes \alpha \in \mathbf{T}$, we have $gx
= x$ iff $x \subset \mathrm{Ker}(\alpha)$ or $\omega \in x$.
\end{lemma}

\dem If $x \subset \mathrm{Ker}(\alpha)$, it is clear that $gx =
x$. Besides, if $x$ is not a subspace of $\mathrm{Ker}(\alpha)$
but $\omega \in x$, we have $x = x' \oplus \mathbb{K} \omega'$
with $x' = x \cap \mathrm{Ker}(\alpha)$ and $\omega' \in x
\setminus \mathrm{Ker}(\alpha)$. This gives $$gx = x' \oplus
\mathbb{K} g (\omega') = x' \oplus \mathbb{K} (\omega' +
\alpha(\omega') \omega) = x' \oplus \mathbb{K} \omega' = x.$$
Reciprocally, if $gx = x$ and $x$ does not belong to
$\mathrm{Ker}(\alpha)$, we take $\omega'' \in x \setminus
\mathrm{Ker}(\alpha)$ so that $g (\omega'') = \omega'' +
\alpha(\omega'') \omega \in x$. Therefore, we conclude that
$\omega \in x$. \fin

Now we combine the expression for the valence of $\X_r$ given in
(\ref{Equation-Valence}) with (\ref{Equation-Cardinal-Graph}) and
Lemma \ref{Lemma-Isotropy} to obtain the value of
$\gamma_1(\mathbf{T}, \X_r)$
\begin{eqnarray*}
\mathrm{val}(\X_r) \, \gamma_1(\mathbf{T}, \X_r) & = & \Big|
\Big\{ g \in \mathbf{T} \, \big| \ gx \neq x \Big\} \Big| \\ & = &
\Big| \Big\{ 1 + \omega \otimes \alpha \, \big| \ \alpha \in
\Omega^* \setminus \{0\}, \omega \in \mathrm{Ker}(\alpha)
\setminus x, x \nsubseteq \mathrm{Ker}(\alpha) \Big\} \Big| \\ & =
& \Big( \Big[ \!\!
\begin{array}{c} n \\ 1 \end{array} \!\! \Big]_q - \Big[ \!\!
\begin{array}{c} n-r \\ 1 \end{array} \!\! \Big]_q \Big) (q^{n-1}
- q^{r-1}).
\end{eqnarray*}
Dividing on the right hand side by the value for
$\mathrm{val}(\X_r)$ given in (\ref{Equation-Valence}), we obtain
the identity $\gamma_1(\mathbf{T}, \X_r) = q^{n-2} (q-1)$. This
proves the first assertion of $\mathrm{(a)}$ in Theorem
\ref{Theorem-Laplacian-Grassmann}. The second assertion follows
from Lemma \ref{Lemma-Laplacian-Ratio}. On the other hand, since
any $\Lambda \in \mathrm{Hom_{\G}}(\V_{r_1},\V_{r_2})$ commutes
with the action of $\G$, we have $$\Lap_{\mathbf{T},{r_2}} \circ
\Lambda = \Lambda \circ \Lap_{\mathbf{T},r_1}.$$ Therefore,
$\mathrm{(b)}$ in Theorem \ref{Theorem-Laplacian-Grassmann}
follows from $\mathrm{(a)}$. Moreover, $\mathrm{(c)}$ is a
consequence of Lemma \ref{Lemma-Laplacian-Irreducible} and
$\mathrm{(a)}$. To prove $\mathrm{(d)}$, we take $(x_1,x_2) \in
\X_{r_1} \times \X_{r_2}$ with $\partial(x_2,x_1) = t$. Then we
have
\begin{eqnarray*}
\Lambda' \varphi (x_2) & = & \sum_{y \in \mathbf{S}_1(x_2)} \big(
\Lambda \varphi(y) - \Lambda \varphi(x_2) \big) \\ & = & \sum_{y
\in \mathbf{S}_1(x_2)} \sum_{x_1 \in \X_{r_1}} \big(
\lambda(\partial(y,x_1)) - \lambda(\partial(x_2,x_1)) \big)
\varphi(x_1).
\end{eqnarray*}
In particular, we can write $$\lambda'(\partial(x_2,x_1)) =
\sum_{y \in \mathbf{S}_1(x_2)} \big( \lambda(\partial(y,x_1)) -
\lambda(\partial(x_2,x_1)) \big).$$ Then $\mathrm{(d)}$ follows
immediately from this. Finally, it remains to see $\mathrm{(e)}$.
The proof requires two combinatorial lemmas. Let us notice that,
given $(x_1,x_2) \in \X_{r_1} \times \X_{r_2}$ and $y \in
\mathbf{S}_1(x_2)$, we have
\begin{eqnarray} \label{Equation-Partial-Dimensions}
\partial(x_2,x_1) - \partial(y,x_1) & = & \partial((x_1+x_2)/x_1) -
\partial((x_1+y)/x_1) \\ \nonumber & = & \partial((x_1+x_2)/(x_1+z)) -
\partial((x_1+y)/(x_1+z)),
\end{eqnarray}
with $z = x_2 \cap y$. Besides, recalling that $y \in
\mathbf{S}_1(x_2)$, it is not difficult to check that both
dimensions appearing on the right hand side of
(\ref{Equation-Partial-Dimensions}) are either $0$ or $1$. This
remark will be used in the following results.

\begin{lemma} \label{Lemma-LinearAlgebra1}
Given $(x_1,x_2) \in \X_{r_1} \times \X_{r_2}$, $(z,\omega) \in
\widetilde{\mathbf{S}}_1(x_2)$ and $y = z \oplus \mathbb{K} \omega
\in \mathbf{S}_1(x_2)$, the following assertions are equivalent:
\begin{itemize}
\item[$\mathrm{(a)}$] $\partial(y,x_1) = \partial(x_2,x_1) + 1.$
\item[$\mathrm{(b)}$] $x_1 + z = x_1 + x_2$ and $\omega \notin
x_1+x_2$.
\item[$\mathrm{(c)}$] $x_1 \cap x_2 \nsubseteq z$ and $\omega \notin
x_1+x_2$.
\end{itemize}
\end{lemma}

\dem Following (\ref{Equation-Partial-Dimensions}) and the remark
after it, we deduce $\mathrm{(a)}$ is equivalent to
$\mathrm{(b)}$. On the other hand, since
$$\partial((x_1+x_2)/(x_1+z)) = 1 - \partial(x_1 \cap x_2) +
\partial(x_1 \cap z),$$ it follows that $\mathrm{(b)}$ is
equivalent to $\mathrm{(c)}$. Therefore, the proof is completed.
\fin

Notice that $z$ is a $(r_2-1)$-dimensional subspace of $x_2$.
Besides, given $x_2 \in \X_{r_2}$, the quotient mapping $\pi: x_2
\rightarrow x_2 / (x_1 \cap x_2)$ provides the following identity
$$\Big| \Big\{ z \subset x_2 \, \big| \ \partial(z) = r_2 -1, x_1
\cap x_2 \subset z \Big\} \Big| = \Big| \Big\{ z' \subset x_2 /
(x_1 \cap x_2) \, \big| \ \partial(z') = \partial(x_2,x_1) - 1
\Big\} \Big|.$$ In other words, the number of subspaces $z$ of
$x_2$ satisfying $\mathrm{(c)}$ is $$\Big[ \!\!
\begin{array}{c} r_2 \\ r_2- 1 \end{array} \!\! \Big]_q - \Big[
\!\! \begin{array}{c} \partial(x_2,x_1) \\
\partial(x_2,x_1) - 1 \end{array} \!\! \Big]_q = \Big[
\!\! \begin{array}{c} r_2 \\ 1 \end{array} \!\! \Big]_q - \Big[
\!\! \begin{array}{c} \partial(x_2,x_1) \\ 1 \end{array} \!\!
\Big]_q.$$ Combining this with Lemmas \ref{Lemma-Covering} and
\ref{Lemma-LinearAlgebra1}, we easily get the following expression
for $\mathbf{b}_{r_1,r_2}$, which simplifies the one given in
Theorem \ref{Theorem-Laplacian-Grassmann}
$$\mathbf{b}_{r_1,r_2}(t) = \Big[ \!\! \begin{array}{c} r_2 \\ 1
\end{array} \!\! \Big]_q \, \frac{q^n - q^{r_1+t}}{q^{r_2-1}(q-1)}
- \Big[ \!\! \begin{array}{c} t \\ 1 \end{array} \!\! \Big]_q
\frac{q^n - q^{r_1+t}}{q^{r_2-1}(q-1)}.$$

\begin{lemma} \label{Lemma-LinearAlgebra2}
Given $(x_1,x_2) \in \X_{r_1} \times \X_{r_2}$, $(z,\omega) \in
\widetilde{\mathbf{S}}_1(x_2)$ and $y = z \oplus \mathbb{K} \omega
\in \mathbf{S}_1(x_2)$, the following assertions are equivalent:
\begin{itemize}
\item[$\mathrm{(a)}$] $\partial(y,x_1) = \partial(x_2,x_1) - 1.$
\item[$\mathrm{(b)}$] $x_1 + z \neq x_1 + x_2$ and $\omega \in
x_1 + z$.
\item[$\mathrm{(c)}$] $x_1 \cap x_2 \subset z$ and $\omega \in
x_1 + z$.
\end{itemize}
Moreover, when them hold we have $\partial(x_1+z)=
\partial(x_1+x_2)-1$ and $z = (x_1 + z) \cap x_2$.
\end{lemma}

\dem By (\ref{Equation-Partial-Dimensions}) and the remark after
it we deduce the equivalence between $\mathrm{(a)}$ and
$\mathrm{(b)}$. The equivalence between $\mathrm{(b)}$ and
$\mathrm{(c)}$ follows again from the identity
$$\partial((x_1+x_2)/(x_1+z)) = 1 - \partial(x_1 \cap x_2) +
\partial(x_1 \cap z).$$ Relation $\partial(x_1+z)=
\partial(x_1+x_2)-1$ is immediate from $\mathrm{(b)}$. The last
claim follows from the modular law $(x_1 + z) \cap x_2 = (x_1 + z)
\cap (x_2 + z) = (x_1 \cap x_2) + z = z$. \fin

Finally, arguing as we did after Lemma \ref{Lemma-LinearAlgebra1},
we easily obtain from Lemmas \ref{Lemma-Covering} and
\ref{Lemma-LinearAlgebra2} the following identity
$$\mathbf{c}_{r_1,r_2}(t) = \Big[ \!\! \begin{array}{c} t \\ 1
\end{array} \!\! \Big]_q \frac{q^{r_1+t-1} -
q^{r_2-1}}{q^{r_2-1}(q-1)}.$$ This is the identity given in
$\mathrm{(e)}$ and the proof of Theorem
\ref{Theorem-Laplacian-Grassmann} is completed.

\section{The kernels in $\mathrm{Hom_G} (\mathrm{V}_{r_1,s},
\mathrm{V}_{r_2,s})$}
\label{Section3}

In this section we obtain explicit formulas for the kernels of the
operators $\Lambda_s^{r_1,r_2}$ in terms of the basic
hypergeometric function. More concretely, it turns out that these
kernels (when regarded as functions of $q^{-t}$ with $t \in
\mathrm{I}_n(r_1,r_2)$) are given by the so-called $q$-Hahn
polynomials. After that, we shall also provide Rodrigues type
formulas for these kernels adapting the techniques developed in
\cite{MP2}. The basic idea consist in showing that the operator
$\Lap_{r_1,r_2}$ introduced above can be identified with the
hypergeometric operator studied in \cite{MP2}.

\subsection{Preliminaries}
\label{Subsection3.1}

As we have pointed out, we shall need some results from the theory
of basic hypergeometric polynomials appearing in \cite{MP2}.
However, since the paper \cite{MP2} considers a great variety of
hypergeometric type operators, we summarize and re-state here
those results from \cite{MP2} that will be useful in the sequel.
To be more precise, we shall formulate the main results from
\cite{MP2} in the particular case of the (so-called there)
\emph{geometric canonical form}. Let us denote by $\Pol_r[x]$ the
space of complex polynomials of degree $\le r$ in one variable.
Given a polynomial $f \in \Pol_r[x]$ and any $q \neq 1$, we define
the linear operators $$\Su f (u) = \frac{f(q^{1/2} u) + f(q^{-1/2}
u)}{2} \qquad \mbox{and} \qquad \Di f(u) = \frac{f(q^{1/2} u) -
f(q^{-1/2} u)}{(q^{1/2} - q^{-1/2}) u}.$$ It is not difficult to
check that $\Su: \Pol_r[x] \rightarrow \Pol_r[x]$ and $\Di:
\Pol_r[x] \rightarrow \Pol_{r-1}[x]$. Moreover, these operators
can be extended so that $\Su, \Di: \mathcal{M}(\C^*) \rightarrow
\mathcal{M}(\C^*)$, where $\mathcal{M}(\C^*)$ stands for the space
of meromorphic functions on $\C \setminus \{0\}$, see
\cite[Section 2.3]{MP2}. Now, given $\sigma \in \Pol_2[x]$ and
$\tau \in \Pol_1[x]$ by $$\sigma(x) = \alpha_2 x^2 + \alpha_1 x +
\alpha_0 \qquad \mbox{and} \qquad \tau(x) = \beta_1 x + \beta_0,$$
we consider the hypergeometric operator $$\OH = \sigma \Di^2 +
\tau \Su \Di.$$ The following result has been adapted from
\cite{MP2} according to our aims.

\begin{lemma} \label{Lemma-TAMS-Hypergeometric}
The hypergeometric operator $\OH: \Pol_r[x] \rightarrow \Pol_r[x]$
satisfy:
\begin{itemize}
\item[$\mathrm{(a)}$] Let us consider two polynomials
\begin{eqnarray*}
\chi^+(u) & = & \gamma_2^+ u^2 + \gamma_1^+ u + \gamma_0, \\
\chi^-(u) & = & \gamma_2^- u^2 + \gamma_1^- u + \gamma_0,
\end{eqnarray*}
with the same value at $u=0$ and let us parameterize $\sigma$ and
$\tau$ as follows $$\begin{array}{rclcrcl} \ \chi^+(u) & = &
\displaystyle \sigma(u) + \frac{q^{-1/2} - q^{1/2}}{2} \, u \,
\tau(u), &  & \sigma(u) & = & \displaystyle \frac{\chi^+(u) +
\chi^-(u)}{2}, \\ \ \chi^-(u) & = & \displaystyle \sigma(u) -
\frac{q^{-1/2} - q^{1/2}}{2} \, u \, \tau(u), &  & \tau(u) & = &
\displaystyle \frac{\chi^+(u) - \chi^-(u)}{(q^{-1/2} - q^{1/2})
u}. \end{array}$$ Assume that the following numbers are pairwise
distinct for $0 \le k \le r$ $$\mu_k = \frac{q^{k/2} -
q^{-k/2}}{(q^{1/2} - q^{-1/2})^2} \left[ \gamma_2^+ q^{(1-k)/2} -
\gamma_2^- q^{(k-1)/2} \right].$$ Then there exist eigenfunctions
$f_k \in \Pol_k[x]$ with degree $k$ satisfying $$\OH f_k(u) +
\mu_k f_k(u) = 0.$$
\item[$\mathrm{(b)}$] Let $f_k$ be as above and let $\partial_k
f_k$ be the main coefficient of $f_k$:

\vskip5pt

\begin{itemize}
\item[$\bullet$] Given $u_0$ such that $\gamma_2^+ u_0^2
+ \gamma_1^+ u_0 + \gamma_0 = 0$, we have $$\quad f_k (u) =
\partial_k f_k \sum_{j=0}^k \Big( \prod_{i=j}^{k-1} \frac{1 -
q^{-(1+i)}}{q^{1/2} (1 - q^{-1})^2} \, \frac{\chi^-(q^i
u_0)}{(\mu_i - \mu_k) u_0} \Big) \prod_{i=0}^{j-1} (u - q^i
u_0).$$
\item[$\bullet$] Given $u_0$ such that $\gamma_2^- u_0^2
+ \gamma_1^- u_0 + \gamma_0 = 0$, we have $$\quad f_k (u) =
\partial_k f_k \sum_{j=0}^k \Big( \prod_{i=j}^{k-1} \frac{1 -
q^{1+i}}{q^{-1/2} (q-1)^2} \, \frac{\chi^+(q^{-i} u_0)}{(\mu_i -
\mu_k) u_0} \Big) \prod_{i=0}^{j-1} (u - q^{-i} u_0).$$
\end{itemize}

\item[$\mathrm{(c)}$] Moreover, assuming above that $f_k(u_0) \neq
0$, we obtain:

\vskip5pt

\begin{itemize}
\item[$\bullet$] If $\gamma_2^+ u_0^2 + \gamma_1^+ u_0 + \gamma_0
= 0$ and $f_k(u_0) \neq 0$, we have $$\qquad f_k (u) = f_k(u_0)
\sum_{j=0}^k \Big( \prod_{i=0}^{j-1} \frac{q^{1/2} (1 -
q^{-1})^2}{1 - q^{-(1+i)}} \, \frac{(\mu_i - \mu_k) u_0}{\chi^-
(q^i u_0)} \Big) \prod_{i=0}^{j-1} (u - q^i u_0).$$
\item[$\bullet$] If $\gamma_2^- u_0^2 + \gamma_1^- u_0 + \gamma_0
= 0$ and $f_k(u_0) \neq 0$, we have $$\qquad f_k(u) = f_k(u_0)
\sum_{j=0}^k \Big( \prod_{i=0}^{j-1} \frac{q^{-1/2} (q-1)^2}{1 -
q^{1+i}} \, \frac{(\mu_i - \mu_k) u_0}{\chi^+ (q^{-i} u_0)} \Big)
\prod_{i=0}^{j-1} (u - q^{-i} u_0).$$
\end{itemize}

\item[$\mathrm{(d)}$] If $\chi^{\pm}(u) = (1 - \xi_0^{\pm} u) (1 -
\xi_1^{\pm} u)$, we obtain the $q$-Hahn polynomials

\vskip5pt

\begin{itemize}
\item[$\bullet$] $\displaystyle f_k(u) = f_k(1/\xi_0^+) \,\ _3\phi_2
\left( \begin{array}{c|} q^k, \ \xi_0^+ u, \ q^{1-k} \xi_0^+
\xi_1^+ / \xi_0^- \xi_1^- \\ \xi_0^+ / \xi_0^-, \ \xi_0^+ /
\xi_1^- \end{array} \ q^{-1}, \ q^{-1} \right).$
\item[$\bullet$] $\displaystyle f_k(u) = f_k(1/\xi_1^+) \,\ _3\phi_2
\left( \begin{array}{c|} q^k, \ \xi_1^+ u, \ q^{1-k} \xi_0^+
\xi_1^+ / \xi_0^- \xi_1^- \\ \xi_1^+ / \xi_0^-, \ \xi_1^+ /
\xi_1^- \end{array} \ q^{-1}, \ q^{-1} \right).$
\item[$\bullet$] $\displaystyle f_k(u) =
f_k(1/\xi_0^-) \,\ _3\phi_2 \left( \begin{array}{c|} q^k, \ 1 /
\xi_0^- u, \ q^{1-k} \xi_0^+ \xi_1^+ / \xi_0^- \xi_1^- \\ \xi_0^+
/ \xi_0^-, \ \xi_1^+ / \xi_0^- \end{array} \ q^{-1}, \ q^{-1}
\xi_1^- u \right)$.
\item[$\bullet$] $\displaystyle f_k(u) =
f_k(1/\xi_1^-) \,\ _3\phi_2 \left( \begin{array}{c|} q^k, \ 1 /
\xi_1^- u, \ q^{1-k} \xi_0^+ \xi_1^+ / \xi_0^- \xi_1^- \\ \xi_0^+
/ \xi_1^-, \ \xi_1^+ / \xi_1^- \end{array} \ q^{-1}, \ q^{-1}
\xi_0^- u \right)$.
\end{itemize}
\end{itemize}
\end{lemma}

\begin{remark}
\emph{The role of $q$ in this paper is played by $q^{-1}$ in
\cite{MP2}.}
\end{remark}

\begin{remark}
\emph{The notation for the basic hypergeometric function follows
\cite{GR}.}
\end{remark}

We shall also give a Rodrigues type formula for the kernels in
$\mathrm{Hom_{\G}}(\V_{r_1,s},\V_{r_2,s})$. To that aim, we state
below the Rodrigues formula given in \cite{MP2} which corresponds
to our problem. That is, the one for the geometric canonical form.

\begin{lemma} \label{Lemma-TAMS-Rodrigues}
Let $\rho \in \mathcal{M}(\C^*)$ be a function satisfying the
functional equation $$\textstyle \rho (u) \chi^+ (u) = \rho(q^{-1}
u) \chi^-(q^{-1} u).$$ Let $\rho_j \in \mathcal{M}(\C^*)$, with $j
\ge 0$ and $\rho_0 = \rho$, determined by any of the recurrences
\begin{eqnarray*}
\rho_{j+1}(u) & = & \textstyle \rho_j(q^{1/2} u) \,
\chi^+(q^{-(j-1)/2} u), \\ \rho_{j+1}(u) & = & \textstyle
\rho_j(q^{-1/2} u) \, \chi^- (q^{(j-1)/2} u).
\end{eqnarray*}
Then, the following Rodrigues formula holds for the eigenfunctions
$f_0, f_1, \ldots, f_r$ $$\Big( \prod_{j=0}^{k-1} (\mu_j - \mu_k)
\Big) \rho(u) f_k (u) = \partial_k f_k \Big( \prod_{j=0}^{k-1}
\frac{q^{(k-j)/2} - q^{(j-k)/2}}{q^{1/2} - q^{-1/2}} \Big) \Di^k
\rho_k(u).$$
\end{lemma}

\subsection{Polynomic expressions}
\label{Subsection3.2}

In this paragraph we express the kernels of the operators of
$\mathbf{B}$ in terms of the basic hypergeometric function. Let us
recall that, given $0 \le r_1, r_2 \le n$, the parameter
$\mathrm{N}(r_1,r_2)$ takes the value $$r_1 \wedge r_2 \wedge
(n-r_1) \wedge (n-r_2)$$ and the dimension of
$\mathrm{Hom_{\G}}(\V_{r_1}, \V_{r_2})$ is $1 +
\mathrm{N}(r_1,r_2)$. In particular, we can define the linear
isomorphism $\Phi: \Pol_{\mathrm{N}(r_1,r_2)} \rightarrow
\mathrm{Hom_{\G}}(\V_{r_1}, \V_{r_2})$ given by $$\Lambda \varphi
(x_2) = \sum_{x_1 \in \X_{r_1}}^{\null} f(q^{-
\partial(x_2,x_1)}) \, \varphi(x_1),$$ for $\Lambda = \Phi f$. In
other words, the kernel $\lambda: \mathrm{I}_n(r_1,r_2)
\rightarrow \C$ of $\Phi f$ has the form $$\lambda(t) =
f(q^{-t}).$$ Then we define $\OH_{r_1,r_2}:
\Pol_{\mathrm{N}(r_1,r_2)} \rightarrow \Pol_{\mathrm{N}(r_1,r_2)}$
by the relation $\OH_{r_1,r_2} = \Phi^{-1} \circ \Lap_{r_1,r_2}
\circ \Phi$.

\begin{lemma} \label{Lemma-Lap-OH}
Let us consider the following polynomials in $\Pol_2[x]$
\begin{eqnarray*}
\chi^+(u) & = & q^{r_1-r_2-1/2} (1 - q^{r_2} u)(1 - q^{n-r_1} u),
\\ \chi^-(u) & = & q^{r_1-r_2-1/2} (1 - u)(1 - q^{r_2-r_1}
u),
\end{eqnarray*}
and let us parameterize $\sigma \in \Pol_2[x]$ and $\tau \in
\Pol_1[x]$ as follows $$\begin{array}{rclcrcl} \chi^+(u) & = &
\displaystyle \sigma(u) + \frac{q^{-1/2} - q^{1/2}}{2} \, u \,
\tau(u), &  & \sigma(u) & = & \displaystyle \frac{\chi^+(u) +
\chi^-(u)}{2}, \\ \chi^-(u) & = & \displaystyle \sigma(u) -
\frac{q^{-1/2} - q^{1/2}}{2} \, u \, \tau(u), &  & \tau(u) & = &
\displaystyle \frac{\chi^+(u) - \chi^-(u)}{(q^{-1/2} - q^{1/2})
u}. \end{array}$$ Then, if $\Su$ and $\Di$ stand for the operators
defined above, we have $$\OH_{r_1,r_2} = \sigma \Di^2 + \tau \Su
\Di.$$
\end{lemma}

\dem Given $f \in \Pol_{\mathrm{N}(r_1,r_2)}$, we have
\begin{eqnarray*}
\lefteqn{\sigma \Di^2 f (q^{-t}) + \tau \Su \Di f (q^{-t})} \\ & =
& \frac{\sigma (q^{-t})}{q^{-t-1/2} - q^{-t+1/2}} \left[
\frac{f(q^{-t-1}) - f(q^{-t})}{q^{-t-1} - q^{-t}} -
\frac{f(q^{-t}) - f(q^{-t+1})}{q^{-t} - q^{-t+1}} \right] \\ & + &
\frac{\tau(q^{-t})}{2} \left[ \frac{f(q^{-t-1}) -
f(q^{-t})}{q^{-t-1} - q^{-t}} + \frac{f(q^{-t}) -
f(q^{-t+1})}{q^{-t} - q^{-t+1}} \right] \\ & = & \frac{2
\sigma(q^{-t}) + (q^{-1/2} - q^{1/2}) q^{-t} \tau(q^{-t})}{2
(q-1)^2 q^{-2t-3/2}} \big( \lambda(t+1) - \lambda(t) \big) \\ & +
& \frac{2 \sigma(q^{-t}) - (q^{-1/2} - q^{1/2}) q^{-t}
\tau(q^{-t})}{2 (q-1)^2 q^{-2t-1/2}} \big( \lambda(t-1) -
\lambda(t) \big) \\ & = & \frac{q^{2t+3/2}
\chi^+(q^{-t})}{(q-1)^2} \big( \lambda(t+1) - \lambda(t) \big) +
\frac{q^{2t+1/2} \chi^-(q^{-t})}{(q-1)^2} \big( \lambda(t-1) -
\lambda(t) \big).
\end{eqnarray*}
Now, applying Theorem \ref{Theorem-Laplacian-Grassmann}, the last
expression equals $$\mathbf{b}_{r_1,r_2}(t) \big( \lambda(t+1) -
\lambda(t) \big) + \mathbf{c}_{r_1,r_2}(t) \big( \lambda(t-1) -
\lambda(t) \big).$$ Applying again Theorem
\ref{Theorem-Laplacian-Grassmann}, the proof is concluded by the
definition of $\OH_{r_1,r_2}$. \fin

\begin{theorem} \label{Theorem-q-Hahn}
The spaces $\V_{r,s}$ are $\Lap_r$-eigenspaces with eigenvalue $-
\mu_s(n)$, where $$\mu_s(n) = (q^s - 1) \frac{q^{n-s+1} -
1}{(q-1)^2}.$$ Moreover, the operators $\Lambda_s^{r_1,r_2} \in
\mathbf{B}$ satisfy $$\Lap_{r_1,r_2} \Lambda_s^{r_1,r_2} +
\mu_s(n) \Lambda_s^{r_1,r_2} = 0.$$ On the other hand, if we
consider the polynomials $f_s^{r_1,r_2} \in \Pol_s[x]$ determined
by $$f_s^{r_1,r_2}(q^{-t}) = \lambda_s^{r_1,r_2}(t),$$ then
$f_s^{r_1,r_2}$ is a $q$-Hahn polynomial of degree $s$ given by
\addtocounter{equation}{-8}
\renewcommand{\theequation}{$\mathsf{H}_q(\arabic{equation})$}
\begin{eqnarray}
& & f_s^{r_1,r_2}(u) = f_s^{r_1,r_2}(1) \,\ _3\phi_2 \left(
\begin{array}{c|} q^s, \ u^{-1}, \ q^{n-s+1} \\ q^{n-r_1}, \
q^{r_2} \end{array} \ q^{-1}, \ q^{r_2-r_1-1}u \right). \\ & &
f_s^{r_1,r_2}(u) = f_s^{r_1,r_2}(q^{-r_2}) \,\ _3\phi_2 \left(
\begin{array}{c|} q^s, \ q^{r_2} u, \ q^{n-s+1} \\ q^{r_1}, \
q^{r_2} \end{array} \ q^{-1}, \ q^{-1} \right). \\ & &
f_s^{r_1,r_2}(u) = f_s^{r_1,r_2}(q^{r_1-n}) \,\ _3\phi_2 \left(
\begin{array}{c|} q^s, \ q^{n-r_1} u, \ q^{n-s+1} \\ q^{n-r_1}, \
q^{n-r_2} \end{array} \ q^{-1}, \ q^{-1} \right). \\ & &
f_s^{r_1,r_2}(u) = f_s^{r_1,r_2}(q^{r_1-r_2}) \,\ _3\phi_2 \left(
\begin{array}{c|} q^s, \ q^{r_1-r_2} u^{-1}, \ q^{n-s+1} \\
q^{r_1}, \ q^{n-r_2} \end{array} \ q^{-1}, \ q^{-1}u \right).
\end{eqnarray}
\end{theorem}
\addtocounter{equation}{4}
\renewcommand{\theequation}{\arabic{equation}}

\dem Notice that $\mu_s(n)$ are pairwise distinct for $0 \le s \le
\mathrm{N}(r_1,r_2)$, since $$\mu_{s_1}(n) - \mu_{s_2}(n) =
q^{s_2} (q^{s_1 - s_2}-1) \frac{q^{n+1 - s_1 - s_2} - 1
}{(q-1)^2}.$$ Therefore, we know from Lemmas
\ref{Lemma-TAMS-Hypergeometric} and \ref{Lemma-Lap-OH} that
$$\left\{ (1 - q^s) \frac{q^{n-s+1} - 1}{(q-1)^2} \, \big| \ 0 \le
s \le \mathrm{N}(r_1,r_2) \right\}$$ is the family of eigenvalues
of $\OH_{r_1,r_2}$. In particular, it turns out that this family
is the family of eigenvalues of $\Lap_r$ when $0 \le s \le r
\wedge (n-r)$ and of $\Lap_{r_1,r_2}$ when $0 \le s \le
\mathrm{N}(r_1,r_2)$. By Theorem \ref{Theorem-Laplacian-Grassmann}
we deduce that, for any $0 \le r \le [n/2] - 1$, all the
eigenvalues of $\Lap_r$ are eigenvalues of $\Lap_{r+1}$ and the
operator $\Lap_{r+1}$ has one more eigenvalue associated to the
eigenspace $\V_{r+1,r+1}$. Applying a simple induction argument,
we know that the eigenvalue $- \mu_s(n)$ of $\Lap_r$ is associated
to the eigenspace $\V_{r,s}$. As a particular case, we obtain the
relation $$\Lap_{r_1,r_2} \Lambda_s^{r_1,r_2} + \mu_s(n)
\Lambda_s^{r_1,r_2} = 0.$$ Once we have identified the eigenvalue
corresponding to the operator $\Lambda_s^{r_1,r_2}$, the given
expressions for the polynomial $f_s^{r_1,r_2}$ in terms of the
$q$-Hahn polynomials follow easily from Lemma
\ref{Lemma-TAMS-Hypergeometric}. This completes the proof. \fin

\begin{remark} \label{Remark-Main-Coefficient}
\emph{Let us denote by $\partial_s f_s^{r_1,r_2}$ the main
coefficient of $f_s^{r_1,r_2 }$. Then, by looking at the main
coefficients of the expressions given in Theorem
\ref{Theorem-q-Hahn}, it is not difficult to check that the
following relations hold
\begin{eqnarray*}
f_s^{r_1,r_2}(1) & = & \partial_s f_s^{r_1,r_2} \, (-1)^s
q^{s(r_1-r_2) - {{s} \choose {2}}} \,
\frac{(q^{n-r_1},q^{r_2};q^{-1})_s}{(q^{n-s+1};q^{-1})_s},
\\ f_s^{r_1,r_2}(q^{-r_2}) & = & \partial_s f_s^{r_1,r_2} \,
q^{-s r_2} \,
\frac{(q^{r_1},q^{r_2};q^{-1})_s}{(q^{n-s+1};q^{-1})_s},
\\ f_s^{r_1,r_2}(q^{r_1-n}) & = & \partial_s f_s^{r_1,r_2} \,
q^{-s(n-r_1)} \,
\frac{(q^{n-r_1},q^{n-r_2};q^{-1})_s}{(q^{n-s+1};q^{-1})_s}. \\
f_s^{r_1,r_2}(q^{r_1-r_2}) & = & \partial_s f_s^{r_1,r_2} \,
(-1)^s q^{- {{s} \choose {2}}} \,
\frac{(q^{r_1},q^{n-r_2};q^{-1})_s}{(q^{n-s+1};q^{-1})_s}.
\end{eqnarray*}}
\end{remark}

\begin{remark}
\emph{The basic hypergeometric series in Theorem
\ref{Theorem-q-Hahn} must be truncated at degree $s$. That is, the
terms of degree $> s$ must be ignored. This is a consequence of
the term $q^s$ which appears in any of them. For instance, we
have}
\begin{eqnarray*}
f_s^{r,s}(u) & = & f_s^{r,s}(1) \,\ _3\phi_2 \left(
\begin{array}{c|} q^s, \ u^{-1}, \ q^{n-s+1} \\ q^{n-r}, \
q^s \end{array} \ q^{-1}, \ q^{s-r-1}u \right) \\ & = &
f_s^{r,s}(1) \, \sum_{k=0}^s \frac{(u^{-1},
q^{n-s+1};q^{-1})_k}{(q^{-1}, q^{n-r}; q^{-1})_k} q^{k(s-r-1)}
u^k,
\end{eqnarray*}
\emph{which is a truncated $_2\phi_1$ series. Now, evaluating at
$u = q^{-t}$ for some $t \in \mathrm{I}_n(r,s)$ and applying the
$q$-Gauss summation formula, we easily obtain} $$f_s^{r,s}(q^{-t})
= f_s^{r,s}(1) \,
\frac{(q^{s-r-1};q^{-1})_t}{(q^{n-r};q^{-1})_t}.$$
\emph{Similarly, we have for $t \in \mathrm{I}_n(n-s,r)$
$$f_s^{n-s,r}(q^{-t}) = f_s^{n-s,r}(1) \,
\frac{(q^{r+s-n-1};q^{-1})_t}{(q^r;q^{-1})_t}.$$ At this point,
the identity $$\frac{f_s^{r_1,r_2}(u)}{f_s^{r_1,r_2}(1)} =
q^{s(r_2-r_1)} \frac{(q^{r_1},q^{n-r_2};
q^{-1})_s}{(q^{n-r_1},q^{r_2};q^{-1})_s} \
\frac{f_s^{r_2,r_1}(q^{r_2-r_1} u)}{f_s^{r_2,r_1}(1)}$$ follows
from certain transformation formulas for the basic hypergeometric
series $_3\phi_2$, see \cite{AAR}. In Section \ref{Section5}, we
shall provide an alternative (combinatorial) proof of this
identity. Therefore, we prefer to omit the details of the proof
just sketched.}
\end{remark}

\subsection{Rodrigues formula}
\label{Subsection3.3}

In this paragraph we provide a Rodrigues type formula for the
kernels $\lambda_s^{r_1,r_2}$ of the operators in $\mathbf{B}$.
More concretely, given any two integers $0 \le r_1, r_2 \le n$ and
$0 \le s \le \mathrm{N}(r_1,r_2)$, we shall study the
eigenfunctions $f_s^{r_1,r_2}$ defined above by the relation
$$f_s^{r_1,r_2}(q^{-t}) = \lambda_s^{r_1,r_2}(t).$$ As it was
noticed in \cite{MP2}, the Rodrigues formula provided by Lemma
\ref{Lemma-TAMS-Rodrigues} is not unique since the given
functional equation has multiple solutions. Hence, the main
difficulty will be to choose the \emph{right} solution of the
functional equation according to our further purposes. Following
Lemma \ref{Lemma-TAMS-Rodrigues}, let us consider $\rho^{r_1,r_2}
\in \mathcal{M}(\C^*)$ satisfying the functional equation
\begin{equation} \label{Equation-FE1}
\rho^{r_1,r_2} (u) \chi^+ (u) = \rho^{r_1,r_2} (q^{-1} u)
\chi^-(q^{-1} u),
\end{equation}
with $\chi^+$ and $\chi^-$ determined by Lemma \ref{Lemma-Lap-OH}.
Moreover, let $$\Big\{ \rho_s^{r_1,r_2} \, \big| \ 0 \le s \le
\mathrm{N}(r_1,r_2) \Big\}$$ be the family of functions in
$\mathcal{M}(\C^*)$ defined by any of the recurrences
\begin{equation} \label{Equation-Recurrence}
\begin{array}{rcl} \rho_{s+1}^{r_1,r_2}(u) & = &
\rho_s^{r_1,r_2}(q^{1/2} u) \, \chi^+(q^{-(s-1)/2} u), \\
\rho_{s+1}^{r_1,r_2}(u) & = & \rho_s^{r_1,r_2}(q^{-1/2} u) \,
\chi^- (q^{(s-1)/2} u),
\end{array}
\end{equation}
where $\rho_0^{r_1,r_2} = \rho^{r_1,r_2}$. Then, implementing in
Lemma \ref{Lemma-TAMS-Rodrigues} the eigenvalues provided by
Theorem \ref{Theorem-q-Hahn}, it is not difficult to see that we
obtain the following Rodrigues formula for the eigenfunctions
$f_s^{r_1,r_2}$
\renewcommand{\theequation}{$\mathsf{R}_s(r_1,r_2)$}
\addtocounter{equation}{-1}
\begin{equation} \label{Equation-Rodrigues}
\rho^{r_1,r_2}(u) f_s^{r_1,r_2}(u) = \partial_s f_s^{r_1,r_2} q^{-
\frac{3}{2} {{s} \choose {2}}} \frac{(q-1)^s}{(q^{n-s+1};
q^{-1})_s} \Di^s \rho_s^{r_1,r_2} (u).
\end{equation}

\renewcommand{\theequation}{\arabic{equation}}

\begin{remark} \label{Remark-Simple-Rodrigues}
\emph{Following Section 4.5 of \cite{MP2}, the easiest solution to
the functional equation (\ref{Equation-FE1}) and the recurrences
(\ref{Equation-Recurrence}) is given by}
\begin{eqnarray*} \varrho^{r_1,r_2}(u) & = & \displaystyle
\frac{(q^{-1}u, q^{r_2-r_1-1}u;q^{-1})_{\infty}}{(q^{r_2}u,
q^{n-r_1}u;q^{-1})_{\infty}}, \\ \varrho_s^{r_1,r_2}(u) & = &
\displaystyle q^{s(r_1-r_2-1/2)} \frac{(q^{-1+s/2}u,
q^{r_2-r_1-1+s/2}u;q^{-1})_{\infty}}{(q^{r_2-s/2}u,
q^{n-r_1-s/2}u;q^{-1})_{\infty}}.
\end{eqnarray*}
\end{remark}

Although the system of functions given in Remark
\ref{Remark-Simple-Rodrigues} provides the simplest Rodrigues
formula for the eigenfunctions $f_s^{r_1,r_2}$, it is not the most
appropriate for our aims. Namely, let us analyze the singularities
of the function $\varrho^{r_1,r_2}$. If we take $u_0 = q^{-t}$ for
some $t \in \Z$, then the function $\varrho^{r_1,r_2}$ has:
\begin{itemize}
\item A double pole at $u_0$ when $$0 \vee (r_2-r_1) \le t \le r_2
\wedge (n-r_1).$$
\item A simple pole at $u_0$ when $$0 \wedge (r_2-r_1) \le t < 0
\vee (r_2-r_1),$$ $$r_2 \wedge (n-r_1) < t \le r_2 \vee (n-r_1).$$
\item A non-vanishing regular point at $u_0$ when $$t < 0 \wedge
(r_2-r_1),$$ $$t > r_2 \vee (n-r_1).$$
\end{itemize}
The function $\varrho^{r_1,r_2}$ does not have any other zeros or
poles in $\C \setminus \{0\}$. In particular, it turns out that
the function $\varrho^{r_1,r_2}$ has singular points at $q^{-t}$
with $t$ belonging to the domain $$\mathrm{I}_n(r_1,r_2) = \Big\{
\partial(x_2,x_1) \, \big| \ x_1 \in \X_{r_1}, x_2 \in \X_{r_2}
\Big\} = \Big\{0 \vee (r_2 - r_1) \le t \le r_2 \wedge (n-r_1)
\Big\}.$$

Obviously, if we want to apply Rodrigues formula
(\ref{Equation-Rodrigues}), we need regular solutions of the
functional equation (\ref{Equation-FE1}) at $q^{-t}$ for $t \in
\mathrm{I}_n(r_1,r_2)$. Any other system of solutions can be
constructed by taking $$\rho^{r_1,r_2}(u) = k(u) \varrho^{r_1,r_2}
(u),$$ with $k$ being meromorphic in $\C \setminus \{0\}$ and
satisfying $k(u) = k(qu)$. Then, the functions $\rho_s^{r_1,r_2}$
arise from $\rho^{r_1,r_2}$ by the recurrences
(\ref{Equation-Recurrence}). The choice of such a function $k$ is
equivalent to the choice of an elliptic function $\mathrm{E}(z) =
k(e^{2 \pi i z})$ with periods $1$ and $$\omega = \frac{1}{2 \pi
i} \log q,$$ see \cite{MP2} for further details. Hence, we need to
find a function $k \in \mathcal{M}(\C^*)$ satisfying $k(u) =
k(qu)$ and having double zeros at $\mathrm{I}_n(r_1,r_2)$. To that
aim, we consider two complex numbers $\xi$ and $\eta$ satisfying
the conditions
\begin{itemize}
\item[i)] The product $\xi \eta$ equals $q^{n-r_1+r_2}$.
\item[ii)] Both $\xi$ and $\eta$ are not of the form $q^{-t}$ for
some integer $t$.
\end{itemize}
Then, the function
\begin{equation} \label{Equation-Elliptic}
k_{r_1,r_2}(u) = \alpha_{r_1,r_2} \frac{(q^{r_2}u,q^{n-r_1}u;
q^{-1})_{\infty}}{(\xi u,\eta u;q^{-1})_{\infty}} \,
\frac{(q^{-r_2-1}u^{-1},q^{-n+r_1-1}u^{-1};
q^{-1})_{\infty}}{(q^{-1}\xi^{-1}u^{-1}, q^{-1}\eta^{-1}u^{-1};
q^{-1})_{\infty}}
\end{equation}
satisfies the required properties. Indeed, the condition $k(u) =
k(qu)$ can be easily checked with the aid of property i). On the
other hand, it follows from property ii) that $k$ has no poles in
$\C \setminus \{0\}$. Therefore, $q^{-t}$ is a double zero of $k$
for any integer $t$ and $k$ has no other zeros. In particular, the
function $$\rho^{r_1,r_2} = k_{r_1,r_2} \, \varrho^{r_1,r_2}$$ is
regular in $\C \setminus \{0\}$, non-vanishing in
$\mathrm{I}_n(r_1,r_2)$ and only vanishes in $q^{-(\Z \setminus
\mathrm{I}_n(r_1,r_2))}$.

\begin{remark} \label{Remark-Alpha(r1,r2)}
\emph{Notice that $\rho^{r_1,r_2}$ is determined up to a constant
$\alpha_{r_1,r_2}$.}
\end{remark}

\begin{theorem} \label{Theorem-Rodrigues}
Given $0 \le r_1,r_2 \le n$, there exists a family of functions
$$\Big\{ \rho_s^{r_1,r_2} \, \big| \ 0 \le s \le
\mathrm{N}(r_1,r_2) \Big\}$$ in $\mathcal{M}(\C^*)$ satisfying the
following properties:
\begin{itemize}
\item[$\mathrm{(a)}$] The function $\rho^{r_1,r_2} =
\rho_0^{r_1,r_2}$ solves the functional equation
$(\ref{Equation-FE1})$.
\item[$\mathrm{(b)}$] The function $\rho^{r_1,r_2}$ is regular
in $\C \setminus \{0\}$ and vanishes in $$q^{-(\Z \setminus
\mathrm{I}_n(r_1,r_2))}.$$
\item[$\mathrm{(c)}$] Given any integer $t \in
\mathrm{I}_n(r_1,r_2)$, we have $$\rho^{r_1,r_2}(q^{-t}) =
q^{t(r_1-r_2+t+1)} \Big[ \!\! \begin{array}{c} n \\
t,r_2-t,r_1-r_2+t,n-r_1-t
\end{array} \!\! \Big]_q.$$
\item[$\mathrm{(d)}$] Each function $\rho_s^{r_1,r_2}$ arise from
$\rho^{r_1,r_2}$ and the recurrences
$(\ref{Equation-Recurrence})$.
\item[$\mathrm{(e)}$] Each function $\rho_s^{r_1,r_2}$ is regular
in $\C \setminus \{0\}$ and vanishes in $$q^{-(\Z \setminus
\mathrm{I}_{n-2s}(r_1-s,r_2-s)) - s/2}.$$
\item[$\mathrm{(f)}$] Given any integer $t \in
\mathrm{I}_{n-2s}(r_1-s,r_2-s)$, we have
$$\rho_s^{r_1,r_2}(q^{-t-s/2}) = \psi_s^{r_1,r_2}(t) \, \Big[ \!\!
\begin{array}{c} n-2s \\ t,r_2-s-t,r_1-r_2+t,n-r_1-s-t \end{array}
\!\! \Big]_q,$$ with $\psi_s^{r_1,r_2}$
given by $$\psi_s^{r_1,r_2}(t) = q^{s(r_1-r_2-1/2) +
t(r_1-r_2+t+1)} (q^n;q^{-1})_{2s}.$$
\end{itemize}
\end{theorem}

\dem Our choice for $\rho^{r_1,r_2}$ will be $k_{r_1,r_2}
\varrho^{r_1,r_2}$ with $k_{r_1,r_2}$ given by
(\ref{Equation-Elliptic}) and the constant $\alpha_{r_1,r_2}$ to
be fixed. Properties (a) and (b) have already been justified. To
prove (c), we observe that the functional equation
(\ref{Equation-FE1}) can be rewritten as
\begin{equation} \label{Equation-FE2}
\rho^{r_1,r_2}(q^{-t}) = \rho^{r_1,r_2}(q^{-t-1})
\frac{\chi^-(q^{-t-1})}{\chi^+(q^{-t})},
\end{equation}
when $t$ and $t+1$ belong to $\mathrm{I}_n(r_1,r_2)$. Therefore,
if we see that the function $$\gamma^{r_1,r_2} (q^{-t}) =
q^{t(r_1-r_2+t+1)} \Big[ \!\! \begin{array}{c} n \\
t,r_2-t,r_1-r_2+t,n-r_1-t
\end{array} \!\! \Big]_q$$ satisfies (\ref{Equation-FE2}), we will
have $\rho^{r_1,r_2}(q^{-t}) = \beta_{r_1,r_2}
\gamma^{r_1,r_2}(q^{-t})$ for any $t \in \mathrm{I}_n(r_1,r_2)$
and some constant $\beta_{r_1,r_2}$. Then, property (c) follows by
taking the appropriate constant $\alpha_{r_1,r_2}$. Let us show
that $\gamma^{r_1,r_2}$  satisfies the functional equation
(\ref{Equation-FE2})
\begin{eqnarray*}
\gamma^{r_1,r_2}(q^{-t}) & = & \frac{q^{t(r_1-r_2+t+1)}
(q;q)_n}{(q;q)_t(q;q)_{r_2-t}(q;q)_{r_1-r_2+t}(q;q)_{n-r_1-t}} \\
& = & \gamma^{r_1,r_2}(q^{-t-1}) q^{-(r_1-r_2+2t+2)}
\frac{(1-q^{t+1})(1-q^{r_1-r_2+t+1})}{(1-q^{r_2-t})(1-q^{n-r_1-t})}
\\ & = & \gamma^{r_1,r_2}(q^{-t-1}) \frac{(1-q^{-t-1})
(1-q^{r_2-r_1-t-1})}{(1-q^{r_2-t})(1-q^{n-r_1-t})} \\ & = &
\gamma^{r_1,r_2}(q^{-t-1})
\frac{\chi^-(q^{-t-1})}{\chi^+(q^{-t})}.
\end{eqnarray*}
Now the constant $\alpha_{r_1,r_2}$ is already fixed so that the
functions $\rho_s^{r_1,r_2}$ are completely determined by property
(d). Taking $u = q^{-t-s/2}$ in the first recurrence in
(\ref{Equation-Recurrence}), we obtain the following relation
\begin{equation} \label{Equation-Iteration}
\rho_s^{r_1,r_2}(q^{-t-s/2}) =
\rho_{s-1}^{r_1,r_2}(q^{-t-(s-1)/2}) \chi^+(q^{-s+1-t}).
\end{equation}
Property (e) means that $$\rho_s^{r_1,r_2}(q^{-t-s/2}) = 0 \qquad
\mbox{for} \qquad t \in \Z \setminus
\mathrm{I}_{n-2s}(r_1-s,r_2-s).$$ This follows easily from the
recurrence (\ref{Equation-Iteration}). Indeed, we just need to
observe which are the zeros of $\chi^+$ and apply (b), we leave
the details to the reader. Therefore, it remains to see (f). Let
us consider the functions $$\gamma_s^{r_1,r_2}(q^{-t-s/2}) =
\psi_s^{r_1,r_2}(t) \, \Big[ \!\!
\begin{array}{c} n-2s \\ t,r_2-s-t,r_1-r_2+t,n-r_1-s-t \end{array}
\!\! \Big]_q.$$ It is not difficult to check that each function
$\gamma_s^{r_1,r_2}$ arise from $\gamma_{s-1}^{r_1,r_2}$ and
(\ref{Equation-Iteration}). Therefore, property (f) follows from
(c) and a simple induction argument. \fin

\begin{remark} \label{Remark-Combinatorial-Meaning}
\emph{In Theorem \ref{Theorem-Rodrigues}, we have chosen the
appropriate solutions of the functional equation
(\ref{Equation-FE1}) and the recurrences
(\ref{Equation-Recurrence}) for our further purposes. This can be
justified by the following combinatorial meaning of these
functions. First, by identity (\ref{Equation-Distances}), we have}
$$\rho^{r_1,r_2}(q^{-t}) q^{-t} = \Big| \Big\{ (x_1,x_2) \in
\X_{r_1} \times \X_{r_2} \, \big| \ \partial(x_2,x_1) = t \Big\}
\Big|.$$ \emph{In particular, (\ref{Equation-Cardinal-Graph})
gives}
\begin{equation} \label{Equation-CM1}
\sum_{t \in \Z} \rho^{r_1,r_2}(q^{-t}) q^{-t} = \Big[ \!\!
\begin{array}{c} n \\ r_1
\end{array} \!\! \Big]_q \Big[ \!\! \begin{array}{c} n \\
r_2 \end{array} \!\! \Big]_q.
\end{equation}
\emph{Second, if $\X_{r}^{n-2s}$ denotes the set of
$r$-dimensional of an $(n-2s)$-dimensional vector space $\Omega_s$
over $\mathbb{K}$, we also have} $$\rho_s^{r_1,r_2}(q^{-t-s/2})
q^{-t} = \delta_s^{r_1,r_2} \Big| \Big\{ (x_1,x_2) \in
\X_{r_1-s}^{n-2s} \times \X_{r_2-s}^{n-2s} \, \big| \
\partial(x_2,x_1) = t \Big\} \Big|.$$
\emph{with $\delta_s^{r_1,r_2} = q^{s(r_1-r_2-1/2)}
(q^n;q^{-1})_{2s}$. In particular,}
\begin{eqnarray} \label{Equation-CM2}
\sum_{t \in \Z + \frac{s}{2}} \rho_s^{r_1,r_2}(q^{-t}) q^{-t} & =
& q^{-s/2} \sum_{t \in \Z} \rho_s^{r_1,r_2}(q^{-t-s/2}) q^{-t} \\
\nonumber & = & q^{s(r_1-r_2-1)} (q^n;q^{-1})_{2s} \Big[ \!\!
\begin{array}{c} n - 2s \\ r_1 - s \end{array} \!\! \Big]_q
\Big[ \!\! \begin{array}{c} n - 2s \\ r_2 - s \end{array} \!\!
\Big]_q \\ \nonumber & = &
\frac{(q^{r_1},q^{n-r_1},q^{r_2},q^{n-r_2};q^{-1})_s}{q^{s(r_2-r_1+1)}
(q^n;q^{-1})_{2s}} \Big[ \!\! \begin{array}{c} n \\ r_1
\end{array} \!\! \Big]_q \Big[ \!\! \begin{array}{c} n \\ r_2
\end{array} \!\! \Big]_q.
\end{eqnarray}
\end{remark}

\begin{remark} \label{Remark-Orthogonality-Relations}
\emph{Our choice in Theorem \ref{Theorem-Rodrigues} has also the
following interpretation. The operator $\Lap_{r_1,r_2}$ is clearly
self-adjoint with respect to the Hilbert-Schmidt inner product on
$\mathrm{Hom}_{\G}(\mathrm{V}_{r_1},\mathrm{V}_{r_2})$. In
particular, this property can be rewritten in terms of the
operator $\OH_{r_1,r_2}$ via the mapping $\Phi:
\Pol_{\mathrm{N}(r_1,r_2)} \rightarrow
\mathrm{Hom}_{\G}(\mathrm{V}_{r_1},\mathrm{V}_{r_2})$. Then, it
can be checked that the hypergeometric operator $\OH_{r_1,r_2}$
becomes self-adjoint with respect to the inner product} $$\langle
f,g \rangle = \sum_{t \in \mathrm{I}_n(r_1,r_2)}
\rho^{r_1,r_2}(q^{-t}) q^{-t} f(q^{-t}) \overline{g(q^{-t})}.$$
\emph{The role of the factor $q^{-t}$ in this expression will
become clear in Lemma \ref{Lemma-Parts} below. The reader is
referred to Section 3.3 of \cite{MP2} for more on these
orthogonality relations.}
\end{remark}

\section{The product formula}
\label{Section4}

In this section we study a product formula for the operators
$\Lambda_s^{r_1,r_2}$ in $\mathbf{B}$. To that aim, our first task
is to normalize these operators since there are only determined up
to a constant factor. We choose the normalization provided by
$$f_s^{r_1,r_2}(1) = \lambda_s^{r_1,r_2}(0) = 1.$$ Following
Remark \ref{Remark-Main-Coefficient}, we have
\begin{equation} \label{Equation-Main-Coefficient}
\partial_s f_s^{r_1,r_2} = (-1)^s q^{s(r_2-r_1)+ {{s} \choose
{2}}} \frac{(q^{n-s+1};q^{-1})_s}{(q^{n-r_1},q^{r_2};q^{-1})_s}.
\end{equation}
In particular, Rodrigues formula (\ref{Equation-Rodrigues})
becomes
\begin{equation} \label{Equation-Normalized-Rodrigues}
\rho^{r_1,r_2}(u) f_s^{r_1,r_2}(u) = (-1)^s q^{s(r_2-r_1) -
\frac{1}{2} {{s} \choose {2}}}
\frac{(q-1)^s}{(q^{n-r_1},q^{r_2};q^{-1})_s} \Di^s
\rho_s^{r_1,r_2} (u).
\end{equation}
Since the representations of $\G$ into the spaces
$\mathrm{V}_{r_1,s_1}$ and $\mathrm{V}_{r_2,s_2}$ are equivalent
if and only if $s_1 = s_2$, we clearly have
$$\Lambda_{s_2}^{r_3,r_4} \circ \Lambda_{s_1}^{r_1,r_2} = 0$$
unless $s_1 = s_2$ and $r_2 = r_3$. In particular, in order to
give an explicit formula for the product of two operators in
$\mbox{End}_{\G}(\mathrm{V})$, it suffices to study the products
$\Lambda_s^{r_2,r_3} \circ \Lambda_s^{r_1,r_2}$. In the following
result, we assume by convention that $$\Big[ \!\! \begin{array}{c}
n \\ -1 \end{array} \!\! \Big]_q = 0$$

\begin{theorem} \label{Theorem-Product}
If $0 \le r_1,r_2,r_3 \le n$ and $0 \le s \le \mathrm{N}(r_1,r_2)
\wedge \mathrm{N}(r_2,r_3)$, we have
\renewcommand{\theequation}{$\mathsf{P}_s(r_1,r_2,r_3)$}
\addtocounter{equation}{-1}
\begin{equation} \label{Equation-Product}
\Lambda_s^{r_2,r_3} \circ \Lambda_s^{r_1,r_2} =
\frac{\mbox{\footnotesize $\Big[$ \hskip-6pt
\begin{tabular}{c} $n$
\\ $r_2$ \end{tabular} \hskip-6pt $\Big]_q$}}{\mbox{\footnotesize
$\Big[$ \hskip-6pt \begin{tabular}{c} $n$ \\ $s$ \end{tabular}
\hskip-6pt $\Big]_q$} - \mbox{\footnotesize $\Big[$ \hskip-6pt
\begin{tabular}{c} $n$ \\ $s-1$ \end{tabular} \hskip-6pt $\Big]_q$}}
\Lambda_s^{r_1,r_3}.
\end{equation}
\renewcommand{\theequation}{\arabic{equation}}
\end{theorem}

The proof of Theorem \ref{Theorem-Product} requires two
preliminary lemmas. In the first one we reduce the proof of the
formula (\ref{Equation-Product}) to two particular cases.

\begin{lemma} \label{Lemma-Particular-Cases}
The product formula $($\ref{Equation-Product}$)$ is implied by:
\begin{itemize}
\item $($\ref{Equation-Product}$)$ for $r_1 = r_3$.
\item $($\ref{Equation-Product}$)$ for $r_1 \le r_2 \le r_3$.
\end{itemize}
\end{lemma}

\dem In what follows we shall write $$\mathsf{k}(r,s) =
\frac{\mbox{\footnotesize $\Big[$ \hskip-6pt
\begin{tabular}{c} $n$ \\ $r$ \end{tabular} \hskip-6pt
$\Big]_q$}}{\mbox{\footnotesize $\Big[$ \hskip-6pt
\begin{tabular}{c} $n$ \\ $s$ \end{tabular} \hskip-6pt $\Big]_q$}
- \mbox{\footnotesize $\Big[$ \hskip-6pt
\begin{tabular}{c} $n$ \\ $s-1$ \end{tabular}
\hskip-6pt $\Big]_q$}}.$$ Let $\mathrm{P}(r,s): \mathrm{V}_r
\rightarrow \mathrm{V}_{r,s}$ be the orthogonal projection. Then,
(\ref{Equation-Cardinal-Graph}) and (\ref{Equation-Dimension})
give
\begin{eqnarray*}
\mbox{tr} (\mathrm{P}(r,s)) & = & \dim \mathrm{V}_{r,s} = \Big[
\!\! \begin{array}{c} n \\ s \end{array} \!\! \Big]_q - \Big[ \!\!
\begin{array}{c} n \\ s-1 \end{array} \!\! \Big]_q, \\
\mbox{tr} (\Lambda_s^{r,r}) & = & \sum_{x \in \X_r}
\lambda_s^{r,r}(\partial (x,x)) =  \Big[ \!\!
\begin{array}{c} n \\ r \end{array} \!\! \Big]_q.
\end{eqnarray*}
In particular, it is clear that $$\mathrm{P}(r,s) =
\frac{\mbox{\footnotesize $\Big[$ \hskip-6pt
\begin{tabular}{c} $n$ \\ $s$ \end{tabular} \hskip-6pt $\Big]_q$}
- \mbox{\footnotesize $\Big[$ \hskip-6pt
\begin{tabular}{c} $n$ \\ $s-1$ \end{tabular}
\hskip-6pt $\Big]_q$}}{\mbox{\footnotesize $\Big[$ \hskip-6pt
\begin{tabular}{c} $n$ \\ $r$ \end{tabular} \hskip-6pt
$\Big]_q$}} \, \Lambda_s^{r,r} = \mathsf{k}(r,s)^{-1}
\Lambda_s^{r,r}.$$ Assuming (\ref{Equation-Product}) for
$r_1=r_3$, we claim that
\begin{equation} \label{Equation-Equivalences}
\mathsf{P}_s(r_1,r_3,r_2) \Leftrightarrow
\mathsf{P}_s(r_1,r_2,r_3) \Leftrightarrow
\mathsf{P}_s(r_2,r_1,r_3).
\end{equation}
For instance, we have
\begin{eqnarray*}
\Lambda_s^{r_3,r_2} \circ \Lambda_s^{r_1,r_3} & = &
\mathsf{k}(r_2,s)^{-1} \Lambda_s^{r_3,r_2} \circ
\Lambda_s^{r_2,r_3} \circ \Lambda_s^{r_1,r_2} \\ & = &
\mathsf{k}(r_3,s) \, \mathsf{k}(r_2,s)^{-1} \, \Lambda_s^{r_2,r_2}
\circ \Lambda_s^{r_1,r_2} \\ & = & \mathsf{k}(r_3,s) \,
\mathrm{P}(r_2,s) \circ \Lambda_s^{r_1,r_2} \\ & = &
\mathsf{k}(r_3,s) \, \Lambda_s^{r_1,r_2}.
\end{eqnarray*}
This proves $\mathsf{P}(r_1,r_2,r_3) \Rightarrow
\mathsf{P}(r_1,r_3,r_2)$ under the assumption of
$\mathsf{P}(r_2,r_3,r_2)$. The other implications in
(\ref{Equation-Equivalences}) can be checked in a similar way. On
the other hand, since that the transpositions $(a,b,c) \mapsto
(a,c,b)$ and $(a,b,c) \mapsto (b,a,c)$ generate all permutations
of $(a,b,c)$, it suffices to prove (\ref{Equation-Product}) in the
particular case $r_1 \le r_2 \le r_3$. Therefore, we need only to
assume the two cases stated above.  \fin

\begin{lemma} \label{Lemma-Parts}
Let $f,g \in \mathcal{M}(\C^*)$ so that $f$ is regular at $q^{-t}$
and $g$ is regular at $q^{-t-1/2}$ for all integer $t$. Assume
also that one of the sets $$\Big\{ t \in \Z \, \big| \ f(q^{-t})
\neq 0 \Big\} \qquad \mbox{or} \qquad \Big\{ t \in \Z \, \big| \
g(q^{-t-1/2}) \neq 0 \Big\}$$ is finite. Then we have the
following summation by parts formula $$\sum_{t \in \Z} f(q^{-t})
\Di g(q^{-t})q^{-t} = - \sum_{t \in \Z + \frac{1}{2}} \Di
f(q^{-t}) g(q^{-t}) q^{-t}.$$
\end{lemma}

\dem We have
\begin{eqnarray*}
\sum_{t \in \Z} f(q^{-t}) \Di g(q^{-t})q^{-t} & = &
\frac{1}{q^{-1/2} - q^{1/2}} \ \sum_{t \in \Z} f(q^{-t}) \big(
g(q^{-t-1/2}) - g(q^{-t+1/2}) \big) \\ & = & \frac{1}{q^{-1/2} -
q^{1/2}} \sum_{t \in \Z+\frac{1}{2}} \big( f(q^{-t+1/2}) -
f(q^{-t-1/2}) \big) g(q^{-t}) \\ & = & - \sum_{t \in \Z +
\frac{1}{2}} \Di f(q^{-t}) g(q^{-t}) q^{-t}. \quad \qquad \qquad
\qquad \qquad \qquad \square
\end{eqnarray*}

\vskip6pt

\noindent \textbf{Proof of (\ref{Equation-Product}) for
$r_1=r_3$.} By Schur lemma, the maps $\Lambda_s^{r_2,r_1} \circ
\Lambda_s^{r_1,r_2}$ and $\mathsf{k}(r_2,s) \Lambda_s^{r_1,r_1}$
are proportional. Hence, it suffices to prove that both operators
have the same trace. Arguing as in Lemma
\ref{Lemma-Particular-Cases}, the operator $\Lambda_s^{r_1,r_1}$
has trace $$\Big[ \!\! \begin{array}{c} n \\ r_1 \end{array} \!\!
\Big]_q.$$ In particular, we need to prove that $$\mbox{tr}
(\Lambda_s^{r_2,r_1} \circ \Lambda_s^{r_1,r_2}) =
\frac{\mbox{\footnotesize $\Big[$ \hskip-6pt
\begin{tabular}{c} $n$ \\ $r_1$ \end{tabular} \hskip-6pt
$\Big]_q$} \mbox{\footnotesize $\Big[$ \hskip-6pt
\begin{tabular}{c} $n$ \\ $r_2$ \end{tabular} \hskip-6pt
$\Big]_q$}}{\mbox{\footnotesize $\Big[$ \hskip-6pt
\begin{tabular}{c} $n$ \\ $s$ \end{tabular} \hskip-6pt $\Big]_q$}
- \mbox{\footnotesize $\Big[$ \hskip-6pt
\begin{tabular}{c} $n$ \\ $s-1$ \end{tabular} \hskip-6pt
$\Big]_q$}}.$$ Let $f_s$ be a polynomial of degree $s$ $$f_s(u) =
\partial_s f_s u^s + \ldots$$ where the dots stand for terms of
lower degree. The action of $\Di$ on the main coefficient of $f_s$
is given by $$\Di f_s (u) = \partial_s f_s \frac{q^{s/2} -
q^{-s/2}}{q^{1/2} - q^{-1/2}} u^{s-1} + \ldots,$$ see Section 2.2
of \cite{MP2}. The iteration of this formula leads to
\begin{eqnarray} \label{Equation-s-Derivative}
\Di^s f_s^{r_1,r_2} & = & \partial_s f_s^{r_1,r_2} (-1)^s
q^{-\frac{1}{2}{{s}\choose{2}}} \frac{(q^s;q^{-1})_s}{(q-1)^s} \\
\nonumber & = & q^{s(r_2-r_1) + \frac{1}{2}{{s}\choose{2}}}
\frac{(q^{n-s+1},q^s;q^{-1})_s}{(q-1)^s
(q^{n-r_1},q^{r_2};q^{-1})_s},
\end{eqnarray}
where the last identity follows from
(\ref{Equation-Main-Coefficient}). Now we are ready to compute the
trace of $\Lambda_s^{r_2,r_1} \circ \Lambda_s^{r_1,r_2}$. We begin
by writing this trace in terms of the eigenfunction
$f_s^{r_1,r_2}$ and the Rodrigues function $\rho^{r_1,r_2}$. To
that aim, we recall the combinatorial meaning of this function,
see Remark \ref{Remark-Combinatorial-Meaning}. Then, we use the
Rodrigues formula given in (\ref{Equation-Normalized-Rodrigues})
after the normalization of the functions $f_s^{r_1,r_2}$.
\begin{eqnarray*}
\mbox{tr} (\Lambda_s^{r_2,r_1} \circ \Lambda_s^{r_1,r_2}) & = &
\sum_{x_1,x_2} \lambda_s^{r_1,r_2}(\partial(x_2,x_1))
\lambda_s^{r_2,r_1}(\partial(x_1,x_2)) \\ & = & \, \sum_{t \in \Z}
\rho^{r_1,r_2}(q^{-t}) q^{-t} f_s^{r_1,r_2}(q^{-t})
f_s^{r_2,r_1}(q^{r_2-r_1-t}) \\ & = & (-1)^s \frac{q^{s(r_2-r_1) -
\frac{1}{2} {{s} \choose {2}}}
(q-1)^2}{(q^{n-r_1},q^{r_2};q^{-1})_s} \\ & \times & \sum_{t \in
\Z} f_s^{r_2,r_1}(q^{r_2-r_1-t}) \Di^s \rho_s^{r_1,r_2}(q^{-t})
q^{-t}.
\end{eqnarray*}
Now recall that $$f_s^{r_2,r_1}(q^{r_2-r_1-t}) - q^{s(r_2-r_1)}
f_s^{r_1,r_2}(q^{-t})$$ is a polynomial of degree less that $s$ in
$u = q^{-t}$. In particular, Lemma \ref{Lemma-Parts} gives
$$\sum_{t \in \Z}^{\null} f_s^{r_2,r_1}(q^{r_2-r_1-t}) \Di^s
\rho_s^{r_1,r_2}(q^{-t}) q^{-t} = q^{s(r_2-r_1)} \sum_{t \in \Z}
f_s^{r_2,r_1}(q^{-t}) \Di^s \rho_s^{r_1,r_2}(q^{-t}) q^{-t}.$$
Then, summation by parts and formulas (\ref{Equation-CM2}) and
(\ref{Equation-s-Derivative}) give
\begin{eqnarray*}
\mbox{tr} (\Lambda_s^{r_2,r_1} \circ \Lambda_s^{r_1,r_2}) & = &
\frac{q^{2s(r_2-r_1) - \frac{1}{2} {{s} \choose {2}}}
(q-1)^2}{(q^{n-r_1},q^{r_2};q^{-1})_s} \, \Di^s f_s^{r_2,r_1}
\sum_{t \in \Z + \frac{s}{2}} \rho_s^{r_1,r_2}(q^{-t}) q^{-t} \\ &
= & \frac{q^{s(r_2-r_1)}
(q^{n-s+1},q^s;q^{-1})_s}{(q^{r_1},q^{n-r_1},q^{r_2},q^{n-r_2};q^{-1})_s}
\sum_{t \in \Z + \frac{s}{2}} \rho_s^{r_1,r_2}(q^{-t}) q^{-t} \\ &
= & q^{-s} \frac{(q^{n-s+1},q^s;q^{-1})_s}{(q^n;q^{-1})_{2s}}
\Big[ \!\! \begin{array}{c} n \\ r_1 \end{array} \!\! \Big]_q
\Big[ \!\! \begin{array}{c} n \\ r_2 \end{array} \!\! \Big]_q \\ &
= & \Big[ \!\! \begin{array}{c} n \\ r_1 \end{array} \!\! \Big]_q
\Big[ \!\! \begin{array}{c} n \\ r_2 \end{array} \!\! \Big]_q
\Big/ \Big( \Big[ \!\! \begin{array}{c} n \\ s \end{array} \!\!
\Big]_q - \Big[ \!\! \begin{array}{c} n \\ s-1 \end{array} \!\!
\Big]_q \Big).
\end{eqnarray*}
Notice that, in the use of (\ref{Equation-s-Derivative}), we
interchange the roles of $r_1$ and $r_2$. \fin

\noindent \textbf{Proof of (\ref{Equation-Product}) for $r_1 \le
r_2 \le r_3$.} Given $0 \le r_1 \le r_2 \le n$, we define the
Radon transform $\mathcal{R}_{\subset}^{r_1,r_2}: \mathrm{V}_{r_1}
\rightarrow \mathrm{V}_{r_2}$ as follows
$$\mathcal{R}_{\subset}^{r_1,r_2} \varphi (x_2) = \sum_{x_1
\subset x_2}^{\null} \varphi(x_1).$$ The kernel of this operator
preserves $\partial$ so that $\mathcal{R}_{\subset}^{r_1,r_2} \in
\mathrm{Hom}_{\G}(\mathrm{V}_{r_1},\mathrm{V}_{r_2})$. On the
other hand, given $(x_1,x_3) \in \X_{r_1} \times \X_{r_2}$,
identity (\ref{Equation-Cardinal-3}) gives $$\Big| \Big\{ x_2 \in
\X_{r_2} \, \big| \ x_1 \subset x_2 \subset x_3 \Big\} \Big| =
\Big[ \!\! \begin{array}{c} r_3 - r_1 \\ r_2 - r_1
\end{array} \!\! \Big]_q = \Big[ \!\! \begin{array}{c} r_3
\\ r_2 \end{array} \!\! \Big]_q \Big[ \!\! \begin{array}{c}
r_2 \\ r_1 \end{array} \!\! \Big]_q \Big/ \Big[ \!\!
\begin{array}{c} r_3 \\ r_1 \end{array} \!\! \Big]_q.$$ In
particular,
\begin{equation} \label{Equation-Radon}
\mathcal{R}_{\subset}^{r_2,r_3} \circ
\mathcal{R}_{\subset}^{r_1,r_2} = \frac{\mbox{\footnotesize
$\Big[$ \hskip-6pt \begin{tabular}{c} $r_3$ \\ $r_2$
\end{tabular} \hskip-6pt $\Big]_q$} \mbox{\footnotesize $\Big[$
\hskip-6pt \begin{tabular}{c} $r_2$ \\ $r_1$ \end{tabular}
\hskip-6pt $\Big]_q$}}{\mbox{\footnotesize $\Big[$ \hskip-6pt
\begin{tabular}{c} $r_3$ \\ $r_1$ \end{tabular} \hskip-6pt
$\Big]_q$}} \, \mathcal{R}_{\subset}^{r_1,r_3}.
\end{equation}
Since $\mathcal{R}_{\subset}^{r_1,r_2} \in
\mathrm{Hom}_{\G}(\mathrm{V}_{r_1},\mathrm{V}_{r_2})$, we
decompose it as
\begin{equation} \label{Equation-Radon-Lambda}
\mathcal{R}_{\subset}^{r_1,r_2} = \sum_{s=0}^{\mathrm{N}(r_1,r_2)}
\mathsf{w_s(r_1,r_2)} \, \Lambda_s^{r_1,r_2}.
\end{equation}
To calculate the coefficients $\mathsf{w_s(r_1,r_2)}$, we observe
from (\ref{Equation-Cardinal-Graph})
$$\mbox{tr}(\Lambda_s^{r_2,r_1} \circ
\mathcal{R}_{\subset}^{r_1,r_2}) = \sum_{(x_1,x_2): \, x_1 \subset
x_2} \lambda_s^{r_2,r_1} (\partial(x_1,x_2)) = \Big[ \!\!
\begin{array}{c} n \\ r_2 \end{array} \!\! \Big]_q \Big[ \!\!
\begin{array}{c} r_2 \\ r_1 \end{array} \!\! \Big]_q.$$ On the
other hand, $$\Lambda_s^{r_2,r_1} \circ
\mathcal{R}_{\subset}^{r_1,r_2} = \mathsf{w}_s(r_1,r_2) \,
\Lambda_s^{r_2,r_1} \circ \Lambda_s^{r_1,r_2}.$$ Therefore,
applying $\mathsf{P}_s(r_1,r_2,r_1)$ we obtain
$$\mbox{tr}(\Lambda_s^{r_2,r_1} \circ
\mathcal{R}_{\subset}^{r_1,r_2}) = \mathsf{w}_s(r_1,r_2) \,
\frac{\mbox{\footnotesize $\Big[$ \hskip-6pt \begin{tabular}{c}
$n$ \\ $r_1$ \end{tabular} \hskip-6pt $\Big]_q$}
\mbox{\footnotesize $\Big[$ \hskip-6pt \begin{tabular}{c} $n$
\\ $r_2$ \end{tabular} \hskip-6pt $\Big]_q$}}{\mbox{\footnotesize
$\Big[$ \hskip-6pt \begin{tabular}{c} $n$ \\ $s$
\end{tabular} \hskip-6pt $\Big]_q$} - \mbox{\footnotesize
$\Big[$ \hskip-6pt \begin{tabular}{c} $n$ \\ $s-1$ \end{tabular}
\hskip-6pt $\Big]_q$}}.$$ This result leads to the exact value of
$\mathsf{w}_s(r_1,r_2)$. Then, we can rewrite
(\ref{Equation-Radon}) using identity
(\ref{Equation-Radon-Lambda}). Although we leave the details to
the reader, it is not difficult to check that this gives
$$\sum_{s=0}^{\mathrm{N}(r_1,r_2)} \frac{\mbox{\footnotesize
$\Big( \Big[$ \hskip-6pt \begin{tabular}{c} $n$ \\ $s$
\end{tabular} \hskip-6pt $\Big]_q$} - \mbox{\footnotesize $\Big[$
\hskip-6pt \begin{tabular}{c} $n$ \\ $s-1$ \end{tabular}
\hskip-6pt $\Big]_q \Big)^2$}}{\mbox{\footnotesize $\Big[$
\hskip-6pt \begin{tabular}{c} $n$ \\ $r_2$ \end{tabular}
\hskip-6pt $\Big]_q$}} \, \Lambda_s^{r_2,r_3} \circ
\Lambda_s^{r_1,r_2} = \sum_{s=0}^{\mathrm{N}(r_1,r_2)} \Big( \Big[
\!\! \begin{array}{c} n \\ s \end{array} \!\! \Big]_q - \Big[ \!\!
\begin{array}{c} n \\ s-1 \end{array} \!\! \Big]_q\Big) \,
\Lambda_s^{r_1,r_3}.$$ Now, since these operators are mutually
orthogonal, we identify coefficients. \fin

\begin{corollary} \label{Corollary-Hilbert-Schmidt}
The Hilbert-Schmidt norm of $\Lambda_s^{r_1,r_2}$ is $$\big\|
\Lambda_s^{r_1,r_2} \big\|_{\mathsf{HS}} = \left[ q^{s(r_2-r_1)}
\frac{(q^{r_1},q^{n-r_2};q^{-1})_s}{(q^{n-r_1},q^{r_2};q^{-1})_s}
\, \frac{\mbox{\footnotesize $\Big[$ \hskip-6pt
\begin{tabular}{c} $n$ \\ $r_1$ \end{tabular} \hskip-6pt $\Big]_q$}
\mbox{\footnotesize $\Big[$ \hskip-6pt \begin{tabular}{c} $n$
\\ $r_2$ \end{tabular} \hskip-6pt
$\Big]_q$}}{\mbox{\footnotesize $\Big[$ \hskip-6pt
\begin{tabular}{c} $n$ \\ $s$ \end{tabular} \hskip-6pt
$\Big]_q$} - \mbox{\footnotesize $\Big[$ \hskip-6pt
\begin{tabular}{c} $n$ \\ $s-1$ \end{tabular} \hskip-6pt
$\Big]_q$}} \right]^{1/2}. $$
\end{corollary}

\dem By Schur lemma, there exists some constant
$\mathsf{c}_s(r_1,r_2)$ with $$(\Lambda_s^{r_1,r_2})^* =
\mathsf{c}_s(r_1,r_2) \, \Lambda_s^{r_2,r_1}.$$ If we write this
relation in terms of the kernels,
$$\lambda_s^{r_1,r_2}(\partial(x_2,x_1)) = \mathsf{c}_s(r_1,r_2)
\, \lambda_s^{r_2,r_1}(\partial(x_1,x_2)) \quad \mbox{for all}
\quad (x_1,x_2) \in \X_{r_1} \times \X_{r_2}.$$ This is equivalent
to $$f_s^{r_1,r_2}(u) = \mathsf{c}_s(r_1,r_2) \, f_s^{r_2,r_1}
(q^{r_2-r_1} u).$$ Identifying the main coefficients via
(\ref{Equation-Main-Coefficient}), we obtain the value of
$\mathsf{c}_s(r_1,r_2)$ $$(\Lambda_s^{r_1,r_2})^* = q^{s(r_2-r_1)}
\frac{(q^{r_1},q^{n-r_2};q^{-1})_s}{(q^{n-r_1},q^{r_2};q^{-1})_s}
\, \Lambda_s^{r_2,r_1}.$$ The given expression for the
Hilbert-Schmidt norm arises from Theorem \ref{Theorem-Product}.
\fin

\begin{remark} \label{Remark-Normalizations}
\emph{Our choice of the basis $\mathbf{B}$ follows from the
condition $$f_s^{r_1,r_2}(1) = \lambda_s^{r_1,r_2}(0) = 1.$$ This
normalization is very natural since in this way the kernels
$\lambda_s^{r,r}$ become the spherical functions associated to the
symmetric space $\X_r$. This will be an essential observation in
Section \ref{Section5}. However, there exist other natural
normalizations for $\Lambda_s^{r_1,r_2}$. For instance, Corollary
\ref{Corollary-Hilbert-Schmidt} provides a normalization for which
the basis $\mathbf{B}$ becomes orthonormal with respect to the
Hilbert-Schmidt inner product. Moreover, combining the results
obtained so far it is not difficult to provide the normalization
for which the basis $\mathbf{B}$ is made up of unitary operators.}
\end{remark}

\section{An alternative proof for the product formula}
\label{Section5}

In this section we provide an alternative proof of
(\ref{Equation-Product}) which does not use any tool from the
theory of classical hypergeometric polynomials. In contrast, the
main tools will be the characterization of spherical functions
given in Theorem \ref{Theorem-Spherical} and the Radon transforms
$$\mathcal{R}_{\subset}^{r_1,r_2}: \mathrm{V}_{r_1} \rightarrow
\mathrm{V}_{r_2}$$ defined for $0 \le r_1 \le r_2 \le n$. Along
the proof, we shall obtain some identities for the kernels
$\lambda_s^{r_1,r_2}$ which might be of independent interest.
Before starting the proof, we recall that Lemma
\ref{Lemma-Particular-Cases} does not use any argument from the
theory of classical hypergeometric polynomials. In particular, we
again reduce the proof of the product formula to the proof of
those particular cases.

\begin{remark}
\emph{Along the proof, we shall assume by convention that}
\begin{itemize}
\item[i)] $\displaystyle \Big[ \!\! \begin{array}{c} n \\ r \end{array}
\!\! \Big]_q = \, 0$ for any integer $r$ not satisfying $0 \le r
\le n$.
\item[ii)] $\Lambda_s^{r_1,r_2} = 0$ for any integer $s$ not satisfying
$0 \le s \le \mathrm{N}(r_1,r_2)$.
\end{itemize}
\end{remark}

\subsection{Combinatorial identities}
\label{Subsection5.1}

Let us consider a subspace $x$ of $\Omega$ of codimension
$\partial(\Omega,x) = t$ and let us fix an integer $0 \le k \le
n$. Our first aim is to calculate the number of $r$-dimensional
subspaces $x_r \in \X_r$ of $\Omega$ satisfying $\partial(x_r,x) =
k$. Clearly, this parameter is invariant under the action of $\G$.
In particular, it depends on the codimension of $x$ but not on $x$
itself. Hence we define $$\mathsf{M}(n,r,t,k) = \Big| \Big\{ x_r
\in \X_r \, \big| \ \partial(x_r,x) = k \Big\} \Big|.$$

\begin{lemma} \label{Lemma-Combinatorics}
We have $$\mathsf{M}(n,r,t,k) = \, q^{k(n-t-r+k)} \Big[ \!\!
\begin{array}{c} t \\ k \end{array} \!\! \Big]_q \Big[ \!\!
\begin{array}{c} n-t \\ r-k \end{array} \!\! \Big]_q.$$
\end{lemma}

\dem Let $$\mathsf{A} = \Big\{ (z,w) \in \X_{r-k} \times \X_k \,
\big| \ z \subset x, w \cap x = \{0\} \Big\}.$$ Then, we compute
the cardinality of $\mathsf{A}$ in two different ways. First we
notice that an element $x_r \in \X_r$ satisfies
$\partial(x_r,x)=k$ if and only if it can be written as $x_r = z
\oplus w$ with $(z,w) \in \mathsf{A}$. We have only one possible
choice for $z = x_r \cap x$ while $w$ is any $k$-dimensional
subspace with $x_r = z \oplus w$. Then, it follows from
(\ref{Equation-Suplement}) that $$|\mathsf{A}| =
\mathsf{M}(n,r,t,k) q^{k(r-k)}.$$ On the other hand, we can count
first how many possible $z$'s can we plug in $\mathsf{A}$ by
applying (\ref{Equation-Cardinal-Graph}) with $(n-t,r-k)$ instead
of $(n,r)$. Then we need to count how many $w$'s can we plug in
$\mathsf{A}$. To that aim we notice that, for any such $w$ the
element $y = x \oplus w$ is an $(n-t+k)$-dimensional subspace
containing $x$. The number of possible $y$'s is given by
(\ref{Equation-Cardinal-3}). Finally, we need to count how many
$w$'s do we have for a fixed $y$, which is given again by
(\ref{Equation-Suplement}). In summary, we find that
$$|\mathsf{A}| = \Big[ \!\! \begin{array}{c} n-t \\ r-k
\end{array} \!\! \Big]_q \Big[ \!\! \begin{array}{c} t \\ k
\end{array} \!\! \Big]_q q^{k(n-t)}.$$ Combining the expressions
obtained so far, we obtain the desired result. \fin

In the following result we use Radon transforms and our formula
for $\mathsf{M}(n,r,t,k)$ to obtain some useful relations between
the kernels $\lambda_s^{r_1,r_2}$ corresponding to a fixed value
of the parameter $s$.

\begin{lemma} \label{Lemma-Fixed-s}
Let $0 \le r_1,r_2,r_3 \le n$ and $0 \le s \le n/2$. Then, there
exist absolute constants $\mathsf{c}_0, \mathsf{c}_1,
\mathsf{c}_2, \mathsf{c}_3$ and $\mathsf{c}_4$ such that for any
$t \in \mathrm{I}_n(r_1,r_3)$, we have:
\begin{itemize}
\item[$\mathrm{(a)}$] $\lambda_s^{r_1,r_3}(t) = \mathsf{c}_0 \,
\lambda_s^{r_3,r_1}(r_1-r_3+t)$.

\vskip3pt

\item[$\mathrm{(b)}$] If $r_2 \le r_3$, \\ $\displaystyle \null
\qquad \quad \mathsf{c}_1 \, \lambda_s^{r_1,r_3}(t) =
\sum_{k}^{\null} \mathsf{M}(r_3,r_2,t,k) \,
\lambda_s^{r_1,r_2}(k)$.
\item[$\mathrm{(c)}$] If $r_1 \le r_2$, \\ $\displaystyle \null
\qquad \quad \mathsf{c}_2 , \lambda_s^{r_1,r_3}(t) =
\sum_{k}^{\null} \mathsf{M}(n-r_1,n-r_2,t,k) \,
\lambda_s^{r_2,r_3}(k)$.
\item[$\mathrm{(d)}$] If $r_2 \le r_1$, \\ $\displaystyle \null
\qquad \quad \mathsf{c}_3 \, \lambda_s^{r_1,r_3}(t) =
\sum_{k}^{\null} \mathsf{M}(r_1,r_2,r_1-r_3+t,r_2-r_3+k) \,
\lambda_s^{r_2,r_3}(k)$.
\item[$\mathrm{(e)}$] If $r_3 \le r_2$, \\ $\displaystyle \null
\qquad \quad \mathsf{c}_4 \, \lambda_s^{r_1,r_3}(t) =
\sum_{k}^{\null} \mathsf{M}(n-r_3,n-r_2,r_1-r_3+t,r_1-r_2+k) \,
\lambda_s^{r_1,r_2}(k)$.
\end{itemize}
\end{lemma}

\dem By Schur lemma, $(\Lambda_s^{r_1,r_3})^* = \mathsf{c}_0 \,
\Lambda_s^{r_3,r_1}$ for some constant $\mathsf{c}_0$ independent
of the variable $t$. Then, (a) follows from the relation between
the corresponding kernels and the identity $$\partial(x_1,x_3) =
r_1-r_3 + \partial(x_3,x_1).$$ To prove (b) we write (again by
Schur lemma) $$\mathsf{c}_1 \, \Lam^{r_1,r_3} = \Rad^{r_2,r_3}
\circ \Lambda_s^{r_1,r_2}.$$ Hence, if $(x_1,x_3) \in \X_{r_1}
\times \X_{r_3}$
\begin{eqnarray*}
\mathsf{c}_1 \, \lam^{r_1,r_3}(\partial(x_3,x_1)) & = & \sum_{x_2
\subset x_3} \lam^{r_1,r_2}(\partial(x_2,x_1)) \\ & = & \ \sum_k
\, \Big| \Big\{ x_2 \in \X_{r_2} \, \big| \ x_2 \subset x_3, \,
\partial(x_2,x_1) = k \Big\} \Big| \, \lam^{r_1,r_2}(k).
\end{eqnarray*}
Then we observe that
\begin{eqnarray*}
\mathsf{M}(r_3,r_2,\partial(x_3,x_1),k) & = & \Big| \Big\{ x_2 \in
\X_{r_2} \, \big| \ x_2 \subset x_3, \, \partial(x_2,x_1 \cap x_3)
= k \Big\} \Big|
\\ & = & \Big| \Big\{ x_2 \in \X_{r_2} \, \big| \ x_2 \subset x_3,
\, \partial(x_2,x_1) = k \Big\} \Big|.
\end{eqnarray*}
To prove (c) we write $\mathsf{c}_2 \, \Lam^{r_1,r_3} =
\Lambda_s^{r_2,r_3} \circ \Rad^{r_1,r_2}$ for some absolute
constant $\mathsf{c}_2$. Proceeding as above, this gives
$$\mathsf{c}_2 \, \lam^{r_1,r_3}(\partial(x_3,x_1)) = \sum_k \,
\Big| \Big\{ x_2 \in \X_{r_2} \, \big| \ x_1 \subset x_2, \,
\partial(x_3,x_2) = k \Big\} \Big| \, \lam^{r_2,r_3}(k).$$ To
calculate the coefficient, we work in the dual space $\Omega^*$.
Let $\mathbb{X}$ stand for the set of linear subspaces of
$\Omega^*$ and $\mathbb{X}_r$ the subset of $r$-dimensional
subspaces of $\Omega^*$. Then, if $x^{\bot} \in \mathbb{X}$
denotes the annihilator of a subspace $x \in \X$, we have
$$\partial(x^{\bot},z^{\bot}) = \partial(x^{\bot}) -
\partial(x^{\bot} \cap z^{\bot}) = \partial(x^{\bot}) -
\partial((x+z)^{\bot}) = \partial(x+z) - \partial(x) =
\partial(z,x).$$ In particular,
\begin{eqnarray*}
\lefteqn{\mathsf{M}(n-r_1,n-r_2,\partial(x_3,x_1),k)} \\ & = &
\Big| \Big\{ x_2^{\bot} \in \mathbb{X}_{n-r_2} \, \big| \
x_2^{\bot} \subset x_1^{\bot}, \, \partial(x_2^{\bot},x_1^{\bot}
\cap x_3^{\bot}) = k \Big\} \Big|
\\ & = & \Big| \Big\{ x_2^{\bot} \in \mathbb{X}_{n-r_2} \, \big| \
x_2^{\bot} \subset x_1^{\bot}, \, \partial(x_2^{\bot},x_3^{\bot})
= k \Big\} \Big|
\\ & = & \Big| \Big\{ x_2 \in \X_{r_2} \, \big| \ x_1 \subset x_2,
\, \partial(x_3,x_2) = k \Big\} \Big|.
\end{eqnarray*}
Finally, (d) follows from (a) and (b) while (e) follows from (a)
and (c). \fin

\begin{theorem} \label{Theorem-Kernel-Relations}
The kernels of the operators in the basis $\mathbf{B}$ satisfy:
\begin{itemize}
\item[$\mathrm{(a)}$] If $0 \le t \le s$ and $s \le r \le n-s$,
\begin{eqnarray*}
\lam^{r,s}(t) & = &
\frac{(q^{s-r-1};q^{-1})_t}{(q^{n-r};q^{-1})_t} \, \lam^{r,s}(0),
\\ \lam^{n-s,r}(t) & = &
\frac{(q^{r+s-n-1};q^{-1})_t}{(q^r;q^{-1})_t} \, \lam^{n-s,r}(0).
\end{eqnarray*}
In particular, we have $\lam^{r,s}(0) \neq 0$ and $\lam^{n-s,r}(0)
\neq 0$ for all $s \le r \le n-s$.
\item[$\mathrm{(b)}$] If $0 \le s \le \mathrm{N}(r_1,r_2)$, there
exists a polynomial $f_s^{r_1,r_2}$ of degree $\le s$ uniquely
determined by the condition $$f_s^{r_1,r_2}(q^{-t}) =
\lam^{r_1,r_2}(t) \qquad \mbox{for} \qquad t \in
\mathrm{I}_n(r_1,r_2).$$ Moreover, the degree of $f_s^{r_1,r_2}$
is $s$ and there are non-zero constants $\mathsf{c}_1$ and
$\mathsf{c}_2$ such that
\begin{eqnarray*}
\qquad f_s^{r_1,r_2}(u) & = & \mathsf{c}_1 \sum_{k=0}^s
\frac{(q^{r_2}u,q^{-1})_{s-k} (u^{-1};q^{-1})_k
u^k}{(q^{s-k};q^{-1})_{s-k} (q^k;q^{-1})_k} \, q^{k(r_2-s+k)}
\lam^{r_1,s}(k), \\ & = & \mathsf{c}_2 \sum_{k=0}^s
\frac{(q^{n-r_1}u,q^{-1})_{s-k} (u^{-1};q^{-1})_k
u^k}{(q^{s-k};q^{-1})_{s-k} (q^k;q^{-1})_k} \, q^{k(n-r_1-s+k)}
\lam^{n-s,r_2}(k).
\end{eqnarray*}
\item[$\mathrm{(c)}$] If $0 \le s \le \mathrm{N}(r_1,r_2)$, there
exists a non-zero absolute constant $\mathsf{c}_0$ such that
$$f_s^{r_1,r_2}(u) = \mathsf{c}_0 \,
f_s^{r_2,r_1}(q^{r_2-r_1}u).$$
\item[$\mathrm{(d)}$] If $0 \le s \le \mathrm{N}(r_1,r_2)$, we
have $f_s^{r_1,r_2}(1) \neq 0$ and
\begin{eqnarray*}
\frac{f_s^{r_1,r_2}(q^{-r_2})}{f_s^{r_1,r_2}(1)} & = &
\frac{\lam^{r_1,s}(s)}{\lam^{r_1,s}(0)} = (-1)^s q^{{{s} \choose
{2}} - r_1s} \frac{(q^{r_1};q^{-1})_s}{(q^{n-r_1};q^{-1})_s}, \\
\frac{f_s^{r_1,r_2}(q^{r_1-n})}{f_s^{r_1,r_2}(1)} & = &
\frac{\lam^{n-s,r_2}(s)}{\lam^{n-s,r_2}(0)} = (-1)^s q^{{{s}
\choose {2}} + s(r_2-n)}
\frac{(q^{n-r_2};q^{-1})_s}{(q^{r_2};q^{-1})_s}.
\end{eqnarray*}
In particular, we have
$$\frac{f_s^{r_1,r_2}(q^{r_1-r_2})}{f_s^{r_1,r_2}(1)} = q^{s(r_2-
r_1)}
\frac{(q^{r_1},q^{n-r_2};q^{-1})_s}{(q^{n-r_1},q^{r_2};q^{-1})_s}.$$
\item[$\mathrm{(e)}$] If $0 \le s \le \mathrm{N}(r_1,r_2)$,
$$(\Lam^{r_1,r_2})^* = \frac{f_s^{r_1,r_2}(1)}{f_s^{r_2,r_1}(1)}
\, q^{s(r_2-r_1)}
\frac{(q^{r_1},q^{n-r_2};q^{-1})_s}{(q^{n-r_1},q^{r_2};q^{-1})_s}
\, \Lam^{r_2,r_1}.$$
\end{itemize}
\end{theorem}

\dem By Lemma \ref{Lemma-Combinatorics} and Lemma
\ref{Lemma-Fixed-s} (e) with $(r_1,r_2,r_3) = (r,s,s-1)$, we have
\begin{eqnarray*}
\lefteqn{\frac{q^{n-r-t}-1}{q-1} \, q^{r-s+t+1} \lam^{r,s}(t+1) +
\frac{q^{r-s+t+1}-1}{q-1} \, \lam^{r,s}(t)} \\ & = & \sum_k \Big[
\!\! \begin{array}{c} r-s+t+1 \\ r-s+k \end{array} \!\! \Big]_q
\Big[ \!\! \begin{array}{c} n-r-t \\ n-r-k \end{array} \!\!
\Big]_q q^{(r-s+k)(k-t)} \lam^{r,s}(k) \\ & = & \sum_k
\mathsf{M}(n-s+1,n-s,r-s+t+1,r-s+k) \, \lam^{r,s}(k) \\ & = &
\mathsf{c}_4 \, \lam^{r,s-1}(t) = 0.
\end{eqnarray*}
The first identity in (a) follows by solving the recurrence
$$\lam^{r,s}(t+1) = \frac{1-q^{s-r-t-1}}{1-q^{n-r-t}} \,
\lam^{r,s}(t).$$ The second identity in (a) follows similarly from
Lemma \ref{Lemma-Combinatorics} and Lemma \ref{Lemma-Fixed-s} (d)
with $(r_1,r_2,r_3) = (n-s+1,n-s,r)$. For the first identity in
(b), we use Lemma \ref{Lemma-Combinatorics} and Lemma
\ref{Lemma-Fixed-s} (b) with $(r_1,s,r_2)$ instead of
$(r_1,r_2,r_3)$
\begin{eqnarray*}
\mathsf{c}_1 \, \lam^{r_1,r_2}(t) & = & \sum_k
\mathsf{M}(r_2,s,t,k) \, \lam^{r_1,s}(k) \\ & = & \sum_k \Big[
\!\! \begin{array}{c} t \\ k \end{array} \!\! \Big]_q \Big[ \!\!
\begin{array}{c} r_2-t \\ s-k \end{array} \!\! \Big]_q
q^{k(r_2-t-s+k)} \lam^{r_1,s}(k) \\ & = & \sum_k
\frac{(q^{r_2-t};q^{-1})_{s-k} (q^t;q^{-1})_k
q^{-tk}}{(q^{s-k};q^{-1})_{s-k} (q^k;q^{-1})_k} \, q^{k(r_2-s+k)}
\lam^{r_1,s}(k).
\end{eqnarray*}
Clearly, the right hand side is a polynomial in the variable
$u=q^{-t}$ of degree $\le s$. This gives the first identity in (b)
and proves the existence of such a polynomial.  Uniqueness follows
from $|\mathrm{I}_n(r_1,r_2)| = \mathrm{N}(r_1,r_2) + 1 > s$. The
second identity in (b) follows in a similar way from Lemma
\ref{Lemma-Combinatorics} and Lemma \ref{Lemma-Fixed-s} (c) by
taking $(r_1,n-s,r_2)$ instead of $(r_1,r_2,r_3)$. From (a) and
the first identity in (b), we can write the main coefficient of
$f_s^{r_1,r_2}$ as $$\sum_{k=0}^s \frac{(-1)^{s-k} q^{r_2k -
{{s-k} \choose {2}}}}{(q^{s-k};q^{-1})_{s-k}(q^k;q^{-1})_k} \,
q^{k(r_2-s+k)} \frac{(q^{s-r_1-1};q^{-1})_k}{(q^{n-r_1};q^{-1})_k}
\, \lam^{r_1,s}(0).$$ Therefore, since the coefficients of
$\lam^{r_1,s}(0)$ in this sum are all positive (notice that
$(q^{s-r_1-1};q^{-1})_k$ is positive for all $k \ge 0$ since $s
\le r_1$), we deduce that the main coefficient does not vanish so
that $f_s^{r_1,r_2}$ has degree $s$. This concludes the proof of
(b). Property (c) is another way to write Lemma
\ref{Lemma-Fixed-s} (a). The property $f_s^{r_1,r_2}(1) \neq 0$
follows trivially from (a). The first identity in (d) follows by
evaluating the first identity in (b) at $u=q^{-r_2}$ and $u=1$ and
then applying (a). Similarly, for the second identity in (d), we
evaluate the second identity in (b) at $u=q^{r_1-n}$ and $u=1$
followed by (a). In both identities, the transformation formula
$$(z;q^{-1})_s = (-1)^s q^{-{{s} \choose {2}}} z^s (q^{s-1}
z^{-1};q^{-1})_s$$ is needed. The last identity in (d) follows
from the previous ones and (c). To prove the identity in (e), it
suffices to notice that the constant $\mathsf{c}_0$ in (c) is
given by $f_s^{r_1,r_2}(1) / f_s^{r_2,r_1}(q^{r_2-r_1})$ and apply
(d). This completes the proof. \fin

\subsection{Proof of the cases $r_1 \le r_2 \le r_3$ and $r_1 \ge
r_2 \ge r_3$} \label{Subsection5.2}

We have already seen that $f_s^{r_1,r_2}$ satisfies the condition
$f_s^{r_1,r_2}(1) \neq 0$ for any $0 \le s \le
\mathrm{N}(r_1,r_2)$. Hence, from now on we normalize the
operators $\Lam^{r_1,r_2}$ in $\mathbf{B}$ requiring
$f_s^{r_1,r_2}(1) = 1$. In particular, now Theorem
\ref{Theorem-Kernel-Relations} (e) has the form
\begin{equation} \label{Equation-Adjoint}
(\Lam^{r_1,r_2})^* = \frac{\mathsf{d}(r_1,s)}{\mathsf{d}(r_2,s)}
\, \Lam^{r_2,r_1} \qquad \mbox{with} \qquad \mathsf{d}(r,s) =
q^{sr} \frac{(q^{n-r};q^{-1})_s}{(q^r;q^{-1})_s}.
\end{equation}
Let us consider $0 \le r_1,r_2,r_3 \le n$ with $r_1 \le r_2 \le
r_3$. Then we claim that $$\mathsf{P}(r_3,r_2,r_1) \Rightarrow
\mathsf{P}(r_1,r_2,r_3).$$ Indeed, let us use the same notation as
in Section \ref{Section4} $$\mathsf{k}(r,s) =
\frac{\mbox{\footnotesize $\Big[$ \hskip-6pt
\begin{tabular}{c} $n$ \\ $r$ \end{tabular} \hskip-6pt
$\Big]_q$}}{\mbox{\footnotesize $\Big[$ \hskip-6pt
\begin{tabular}{c} $n$ \\ $s$ \end{tabular} \hskip-6pt $\Big]_q$}
- \mbox{\footnotesize $\Big[$ \hskip-6pt
\begin{tabular}{c} $n$ \\ $s-1$ \end{tabular}
\hskip-6pt $\Big]_q$}}.$$ Then (\ref{Equation-Adjoint}) gives
\begin{eqnarray*}
\Lam^{r_2,r_3} \circ \Lam^{r_1,r_2} & = &
\frac{\mathsf{d}(r_3,s)}{\mathsf{d}(r_2,s)} \, (\Lam^{r_3,r_2})^*
\frac{\mathsf{d}(r_2,s)}{\mathsf{d}(r_1,s)} \, (\Lam^{r_2,r_1})^*
\\ & = & \frac{\mathsf{d}(r_3,s)}{\mathsf{d}(r_1,s)} \,
(\Lam^{r_2,r_1} \circ \Lam^{r_3,r_2})^* \\ & = &
\frac{\mathsf{d}(r_3,s)}{\mathsf{d}(r_1,s)} \, \mathsf{k}(r_2,s)
\, (\Lam^{r_3,r_1})^* = \mathsf{k}(r_2,s) \, \Lam^{r_1,r_3}.
\end{eqnarray*}
In summary, it suffices to prove the case $r_1 \ge r_2 \ge r_3$.
Notice that this is clear since Lemma \ref{Lemma-Particular-Cases}
obviously holds with $r_1 \ge r_2 \ge r_3$ instead of $r_1 \le r_2
\le r_3$. However, we have proved the implication
$\mathsf{P}(r_3,r_2,r_1) \Rightarrow \mathsf{P}(r_1,r_2,r_3)$
since we shall need both results in Paragraph \ref{Subsection5.3}.

\begin{remark} \label{Remark-Spherical}
\emph{Let us write again $\mathrm{P}(r,s): \mathrm{V}_r
\rightarrow \mathrm{V}_{r,s}$ for the orthogonal projection from
$\mathrm{V}_r$ onto $\mathrm{V}_{r,s}$. Then, using that
$\lam^{r,r}(0) = f_s^{r,r}(1) = 1$ and arguing as in the proof of
Lemma \ref{Lemma-Particular-Cases}, we have}
\begin{equation} \label{Equation-Projection-Trace}
\mathrm{P}(r,s) = \frac{\dim \mathrm{V}_{r,s}}{|\X_r|} \,
\Lam^{r,r} = \mathsf{k}(r,s)^{-1} \, \Lam^{r,r}.
\end{equation}
\emph{In particular, it turns out that}
\begin{eqnarray*}
\Lam^{r_1,r_2} \circ \Lam^{r_1,r_1} & = & \mathsf{k}(r_1,s) \,
\Lam^{r_1,r_2}, \\ \Lam^{r_2,r_2} \circ \Lam^{r_1,r_2} & = &
\mathsf{k}(r_2,s) \, \Lam^{r_1,r_2}.
\end{eqnarray*}
\emph{That is, (\ref{Equation-Product}) holds with $r_1 = r_2$ or
$r_2 = r_3$. Moreover, recalling the definition of spherical
function given in Section \ref{Section1} and that $\X_r$ is a
finite symmetric space for any $0 \le r \le n$, we observe from
(\ref{Equation-Projection-Trace}) that the set of spherical
functions associated to $\X_r$ is $$\Big\{ \lam^{r,r} \, \big| \ 0
\le s \le r \wedge (n-r) \Big\}.$$}
\end{remark}

Following Remark \ref{Remark-Spherical}, we are now allowed to use
the characterization of spherical functions provided by Theorem
\ref{Theorem-Spherical}. That is, given $y \in \X_r$ for some $0
\le r \le n$, we consider the isotropy subgroup of $y$ $$\G_y =
\Big\{ g \in \G \, \big| \ gy=y \Big\}.$$ Then, Theorem
\ref{Theorem-Spherical} gives
\begin{equation} \label{Equation-Spherical-Xr}
\frac{1}{|\G_y|} \sum_{g \in \G_y} \lam^{r,r}(\partial(gx,z)) =
\lam^{r,r}(\partial(x,y)) \, \lam^{r,r}(\partial(y,z)),
\end{equation}
for any $x,y,z \in \X_r$. On the other hand, by Remark
\ref{Remark-Spherical} we have
\begin{eqnarray*}
\lam^{r_1,r_2}(\partial(x_2,x_1)) & = &
\frac{1}{\mathsf{k}(r_1,s)} \sum_{z_1 \in \X_{r_1}} \lam^{r_1,r_2}
(\partial(x_2,z_1)) \, \lam^{r_1,r_1} (\partial(z_1,x_1)),
\\ \lam^{r_1,r_2}(\partial(x_2,x_1)) & = &
\frac{1}{\mathsf{k}(r_2,s)} \sum_{z_2 \in \X_{r_2}} \lam^{r_2,r_2}
(\partial(x_2,z_2)) \, \lam^{r_1,r_2} (\partial(z_2,x_1)).
\end{eqnarray*}
Combining these identities with (\ref{Equation-Spherical-Xr}), we
obtain
\begin{eqnarray*}
\frac{1}{|\G_{y_1}|} \sum_{g \in \G_{y_1}}
\lam^{r_1,r_2}(\partial(x_2,gx_1)) & = &
\lam^{r_1,r_2}(\partial(x_2,y_1)) \,
\lam^{r_1,r_1}(\partial(y_1,x_1)), \\ \frac{1}{|\G_{y_2}|} \sum_{g
\in \G_{y_2}} \lam^{r_1,r_2}(\partial(gx_2,x_1)) & = &
\lam^{r_2,r_2}(\partial(x_2,y_2)) \,
\lam^{r_1,r_2}(\partial(y_2,x_1)).
\end{eqnarray*}
with $(y_1,y_2) \in \X_{r_1} \times \X_{r_2}$. Now, let us assume
that $r_1 \ge r_2 \ge r_3$ and let us take $x_k \in \X_{r_k}$ for
$k=1,2,3$ with $x_3 \subset x_2 \subset x_1$. Clearly, for any $g
\in \G_{x_2}$ we will have $x_3 \subset x_2 \subset gx_1$ so that
$\partial(x_3,gx_1)=0$. On the other hand, by Schur lemma we have
$$\Lam^{r_2,r_3} \circ \Lam^{r_1,r_2} = \mathsf{k}_s(r_1,r_2,r_3)
\Lam^{r_1,r_3}$$ for some constant $\mathsf{k}_s(r_1,r_2,r_3)$.
Putting the previous results all together, we have
\begin{eqnarray*}
\mathsf{k}_s(r_1,r_2,r_3) & = &
\frac{\mathsf{k}_s(r_1,r_2,r_3)}{|\G_{x_2}|} \sum_{g \in \G_{x_2}}
\lam^{r_1,r_3}(\partial(x_3,gx_1)) \\ & = & \frac{1}{|\G_{x_2}|}
\sum_{g \in \G_{x_2}} \sum_{z_2 \in \X_{r_2}}
\lam^{r_2,r_3}(\partial(x_3,z_2)) \,
\lam^{r_1,r_2}(\partial(z_2,gx_1)) \\ & = & \frac{1}{|\G_{x_2}|}
\sum_{g \in \G_{x_2}} \sum_{z_2 \in \X_{r_2}}
\lam^{r_2,r_3}(\partial(x_3,z_2)) \,
\lam^{r_1,r_2}(\partial(gz_2,x_1)) \\ & = & \sum_{z_2 \in
\X_{r_2}} \lam^{r_2,r_3}(\partial(x_3,z_2)) \,
\lam^{r_2,r_2}(\partial(z_2,x_2)) \, \lam^{r_1,r_2}
(\partial(x_2,x_1)) \\ & = & \mathsf{k}(r_2,s) \,
\lam^{r_2,r_3}(\partial(x_3,x_2)) \,
\lam^{r_1,r_2}(\partial(x_2,x_1)) = \mathsf{k}(r_2,s).
\end{eqnarray*}

\subsection{Proof of the case $r_1 = r_3$}
\label{Subsection5.3}

In this paragraph we shall need to use another type of Radon
transforms. Given an integer $0 \le r \le n$, we consider the
Radon transform $\mathcal{R}_c^r: \mathrm{V}_r \rightarrow
\mathrm{V}_{n-r}$ defined as follows $$\mathcal{R}_c^r \varphi(z)
= \sum_{x: \ x \cap z = \{0\}} \varphi(x) \qquad \mbox{for} \qquad
z \in \mathrm{V}_{n-r}.$$ Since the kernel of $\mathcal{R}_c^r$ is
invariant under the action of $\G$, $\mathcal{R}_c^r \in
\mathrm{Hom}_{\G}(\mathrm{V}_r,\mathrm{V}_{n-r})$.

\begin{lemma} \label{Lemma-Radon-Complement}
If $r_1 + r_2 \le n$, we have $\mathcal{R}_c^{r_2} \circ
\Lam^{r_1,r_2} = \mathsf{m}(r_2,s) \, \Lam^{r_1,n-r_2}$ with
$$\mathsf{m}(r,s) = (-1)^s q^{(r-s)(n-r) + {{s} \choose {2}}}
\frac{(q^{n-r};q^{-1})_s}{(q^r;q^{-1})_s}.$$
\end{lemma}

\dem By Schur lemma, we know the existence of a constant
$\mathsf{c}$ such that $$\mathcal{R}_c^{r_2} \circ \Lam^{r_1,r_2}
= \mathsf{c} \, \Lam^{r_1,n-r_2}.$$ When $r_1+r_2 \le n$, we can
take $x_1 \in \X_{r_1}$ and $z_2 \in \X_{n-r_2}$ such that $x_1
\subset z_2$. Then, $\partial(z_2,x_1) = n-r_1-r_2$ and identity
(\ref{Equation-Suplement}) gives
\begin{eqnarray*}
\mathsf{c} \, \lam^{r_1,n-r_2} (n-r_1-r_2) & = & \sum_{x_2: \ x_2
\cap z_2 = \{0\}}^{\null} \lam^{r_1,r_2}(\partial(x_2,x_1))
\\ & = & q^{r_2(n-r_2)} \lam^{r_1,r_2}(r_2).
\end{eqnarray*}
Since we have
\begin{eqnarray*}
\lam^{r_1,r_2}(r_2) & = & f_s^{r_1,r_2}(q^{-r_2}), \\
\lam^{r_1,n-r_2} (n-r_1-r_2) & = & f_s^{r_1,n-r_2}(q^{r_1+r_2-n}),
\end{eqnarray*}
the desired relation can be easily checked by applying Theorem
\ref{Theorem-Kernel-Relations} (d). \fin

Now we are ready to complete the proof. Arguing as in Section
\ref{Section4}, the product formula
\begin{equation} \label{Equation-Product-2}
\Lambda_s^{r_2,r_1} \circ \Lambda_s^{r_1,r_2} = \mathsf{k}(r_2,s)
\, \Lam^{r_1,r_1}
\end{equation}
is equivalent to $$\mbox{tr} (\Lambda_s^{r_2,r_1} \circ
\Lambda_s^{r_1,r_2}) = \frac{\mbox{\footnotesize $\Big[$
\hskip-6pt \begin{tabular}{c} $n$ \\ $r_1$ \end{tabular}
\hskip-6pt $\Big]_q$} \mbox{\footnotesize $\Big[$ \hskip-6pt
\begin{tabular}{c} $n$ \\ $r_2$ \end{tabular} \hskip-6pt
$\Big]_q$}}{\mbox{\footnotesize $\Big[$ \hskip-6pt
\begin{tabular}{c} $n$ \\ $s$ \end{tabular} \hskip-6pt $\Big]_q$}
- \mbox{\footnotesize $\Big[$ \hskip-6pt
\begin{tabular}{c} $n$ \\ $s-1$ \end{tabular} \hskip-6pt
$\Big]_q$}}.$$ Since $\mbox{tr} (\Lambda_s^{r_2,r_1} \circ
\Lambda_s^{r_1,r_2}) = \mbox{tr} (\Lambda_s^{r_1,r_2} \circ
\Lambda_s^{r_2,r_1})$, we assume without lost of generality that
$r_1 \le r_2$. In particular, we have $s \le r_1 \le r_2 \le n-s$.
Multiplying on the left (resp. right) of
(\ref{Equation-Product-2}) by $\Lam^{r_2,s}$ (resp.
$\Lam^{s,r_1}$) and applying the results obtained in Paragraph
\ref{Subsection5.2}, it turns out that the proof of
(\ref{Equation-Product-2}) is equivalent to the proof of
$$\Lambda_s^{r_2,s} \circ \Lambda_s^{s,r_2} = \mathsf{k}(r_2,s) \,
\Lam^{s,s}.$$ Now, multiplying on the left by $\mathcal{R}_c^s$
and applying Lemma \ref{Lemma-Radon-Complement} (notice that
$r_2+s \le n$ and $s+s \le n$), the proof of
(\ref{Equation-Product-2}) becomes equivalent to $$\Lam^{r_2,n-s}
\circ \Lam^{s,r_2} = \mathsf{k}(r_2,s) \, \Lam^{s,n-s}.$$ However,
this identity holds since $s \le r_2 \le n-s$. Therefore, the
proof is completed.

\bibliographystyle{amsplain}

\end{document}